\documentclass[11pt,oneside, letter]{article}

\usepackage{amsthm, amsmath, amssymb, amsfonts, graphicx, epsfig}
\usepackage{accents}

\usepackage{bbm}
\usepackage{mathrsfs}

\usepackage{epstopdf}

\usepackage{graphicx}
\usepackage[font=sl,labelfont=bf]{caption}
\usepackage{subcaption}

\usepackage{float}
\restylefloat{table}

\usepackage{enumerate}

\usepackage[colorlinks=true, pdfstartview=FitV, linkcolor=blue,
            citecolor=blue, urlcolor=blue]{hyperref}
\usepackage[usenames]{color}

\usepackage{url}
\makeatletter\def\url@leostyle{%
 \@ifundefined{selectfont}{\def\UrlFont{\sf}}{\def\UrlFont{\scriptsize\ttfamily}}} \makeatother

\usepackage[margin=1.0in, letterpaper]{geometry}

\newtheorem{theorem}{Theorem}
\newtheorem{proposition}[theorem]{Proposition}

\theoremstyle{definition}
\newtheorem{definition}[theorem]{Definition}
\newtheorem{example}[theorem]{Example}
\theoremstyle{remark}
\newtheorem{remark}[theorem]{Remark}

\numberwithin{equation}{section}
\numberwithin{theorem}{section}

\definecolor{Red}{rgb}{0.9,0,0.0}
\definecolor{Blue}{rgb}{0,0.0,1.0}

\def\cB{\mathcal{B}}

\def\cD{\mathcal{D}}

\def\cF{\mathcal{F}}
\def\cG{\mathcal{G}}

\def\cT{\mathcal{T}}

\def\cX{\mathcal{X}}
\def\cY{\mathcal{Y}}

\def\bF{\mathbb{F}}

\def\bN{\mathbb{N}}

\def\bR{\mathbb{R}}

\def\bT{\mathbb{T}}

\def\bV{\mathbb{V}}

\def\bZ{\mathbb{Z}}

\newcommand{\1}{\mathbbm{1}}            
\newcommand{\set}[1]{\{#1\}}            
\newcommand{\Set}[1]{\left\{#1\right\}} 
\renewcommand{\mid}{\;|\;}              

\DeclareMathOperator*{\esssup}{ess\,sup} 
\DeclareMathOperator*{\essinf}{ess\,inf} 

\DeclareMathOperator{\Essinf}{ess\,inf} 
\DeclareMathOperator{\Esssup}{ess\,sup} 

\title{ 
        A unified approach to time consistency of dynamic risk measures and dynamic performance measures in discrete time
}

\makeatletter
\def\and{%
  \end{tabular}%
  \begin{tabular}[t]{c}}%
\def\@fnsymbol#1{\ensuremath{\ifcase#1\or a\or b\or c\or
   d\or e\or f\or g\or h\or i\else\@ctrerr\fi}}
\makeatother

\author{
        Tomasz R. Bielecki\,\thanks{Department of Applied Mathematics, Illinois Institute of Technology
       \newline \hspace*{1.45em}  10 W 32nd Str, Building RE, Room 208, Chicago, IL 60616, USA
       \newline \hspace*{1.45em}  Emails: \url{bielecki@iit.edu} (T.R. Bielecki) and \url{cialenco@iit.edu} (I. Cialenco)
       \newline \hspace*{1.45em}  URLs: \url{http://math.iit.edu/\~bielecki}  and \url{http://math.iit.edu/\~igor}
        \vspace{0.5em}} ,
\and
        Igor Cialenco\,\footnotemark[1] ,
\and and
        Marcin Pitera\,\thanks{Institute of Mathematics, Jagiellonian University,  Cracow, Poland
        \newline \hspace*{1.45em}    Email: \url{marcin.pitera@im.uj.edu.pl}, URL: \url{http://www2.im.uj.edu.pl/MarcinPitera/} }
        }

\date{ {\small
First Circulated: September 22, 2014\\
This Version: January 30, 2017 \\[1em]
Forthcoming in Mathematics of Operations Research \\
\url{http://doi.org/10.1287/moor.2017.0858}
}} %
\begin{document}

\maketitle

\smallskip

{\footnotesize
\begin{tabular}{l@{} p{350pt}}
  \hline \\[-.2em]
  \textsc{Abstract}: \ &
  In this paper we provide a flexible framework allowing for a unified study of time consistency of risk measures and performance measures (also known as acceptability indices). The proposed framework not only integrates existing forms of time consistency, but also provides a comprehensive toolbox for analysis and synthesis of the concept of time consistency in decision making. In particular, it allows for in depth comparative analysis of (most of) the existing types of time consistency -- a feat that has not be possible before and which is done in the companion paper \cite{BCP2014} to this one.
In our approach the time consistency is studied for a large class of maps that are postulated to satisfy only two properties -- monotonicity and locality. The time consistency is defined in terms of an update rule. The form of the update rule introduced here is novel, and is perfectly suited for developing the unifying framework that is worked out in this paper.
As an illustration of the applicability of our approach, we show how to recover almost all concepts of weak time consistency by means of constructing appropriate  update rules. \\[0.5em]
\textsc{Keywords:} \ &  time consistency, update rule, dynamic LM-measure, dynamic risk measure, dynamic acceptability index, dynamic performance measure. \\
\textsc{MSC2010:} \ & 91B30, 62P05, 97M30, 91B06. \\[1em]
  \hline
\end{tabular}
}

\section{Introduction}

In the seminal paper by  Artzner et al. \cite{ArtznerDelbaenEberHeath1999}, the authors proposed an axiomatic approach to  defining risk measures that are meant to give a numerical value of the riskiness of a given financial contract or portfolio. Alternatively, one can view the risk measures as a tool that allows to establish preference orders on the set of cash flows according to their riskiness. Another seminal paper,  Cherny and Madan \cite{ChernyMadan2009}, introduced and studied axiomatic approach to defining performance measures, or acceptability  indices, that are meant to provide evaluation of performance of a financial portfolio. In their most native form, performance measures evaluate the trade-off between return on the portfolio and the portfolio's risk. Both Artzner et al. \cite{ArtznerDelbaenEberHeath1999} and  Cherny and Madan \cite{ChernyMadan2009} were concerned with measures of risk and measures of performance in the static framework.

As shown in one of the first papers that studied  risk measures in the dynamic framework, Riedel \cite{Riedel2004}, if one is concerned about making noncontradictory decisions (from the risk point of view) over  time, then an additional axiom, called time consistency, is needed. Over the past decade significant progress has been made towards expanding the theory of dynamic risk measures and their time consistency. {For example, so called cocycle condition (for convex risk measures) was studied in F\"ollmer and Penner~\cite{FollmerPenner2006}, recursive construction was exploited in Cheridito and Kupper~\cite{CheriditoKupper2006}, relation to acceptance and rejection sets was studied in Delbaen~\cite{Delbaen2006}, the concept of prudence was introduced in Penner~\cite{Penner04}, connections to g-expectations were studied in  Rosazza Gianin~\cite{RosazzaGianin2006}, and the relation to Bellman's principle of optimality was shown in Artzner et al.~\cite{ArtznerDelbaenEberHeathKu2007}.}

{Following Acciaio and Penner \cite{AcciaioPenner2010}} let us briefly recall the concept of strong time consistency of {dynamic} monetary risk measures,\footnote{A dynamic monetary risk measure is a local, monotone and cash-additive function; see Section~\ref{sec:prelims} for a formal definition.} which is one of the most recognized forms of time consistency. Assume that $\rho_t(X)$ is the value of a  monetary risk measure at time $t\in[0,T]$, that corresponds to the riskiness, at time $t$, of the terminal cash flow $X$, with $X$  being an $\cF_T$-measurable random variable. The {dynamic} monetary risk measure {$\rho=\{\rho_t\}_{0\leq t\leq T}$} is said to be strongly time consistent if for any $t<s\leq T$, and any $\cF_T$-measurable random variables $X,Y$ we have that
\begin{equation}\label{eq:strong1}
  \rho_s(X) = \rho_s(Y) \quad \Rightarrow \quad \rho_t(X)= \rho_t(Y).
\end{equation}
The financial interpretation of {the} strong time consistency is clear -- if $X$ is as risky as $Y$ at some future time $s$, then today, at time $t$, $X$ is also as risky as $Y$.
One of the main features of the strong time consistency is its connection to {the} dynamic programming principle. It is not hard to show that in the $L^{\infty}$ framework, a {dynamic} monetary risk measure is strongly time consistent if and only if
\begin{equation}\label{eq:DPP1}
\rho_t = \rho_t (-\rho_s), \quad 0\leq t<s\leq T.
\end{equation}
All other forms of time consistency for {dynamic} monetary risk measures, such as weak, acceptance consistent, rejection consistent,  are tied to this connection as well. In  Tutsch \cite{Tutsch2008}, the author proposed a general approach to time consistency for cash-additive risk measures by introducing so called {\it test sets} or {\it benchmark sets}. Each form of time consistency was associated to a benchmark set of random variables, and larger benchmark sets corresponded to stronger forms of time consistency.

For more details on dynamic cash-additive (monetary risk) measures and their time consistency, we refer the reader to a comprehensive survey paper Acciaio and Penner~\cite{AcciaioPenner2010} and the references therein.

Besides the dynamic risk measures, in this paper we study the dynamic acceptability indices that are also known as dynamic performance measures.\footnote{A dynamic acceptability index is a local, monotone and scale invariant function; see Section~\ref{sec:prelims} for further discussions.}
The scale invariance property, which is the distinctive property of dynamic performance measures, makes the study of time consistency in this case more intricate. In particular, the recursive property analogous to \eqref{eq:DPP1} or the benchmark sets approach  are not appropriate for study of time consistency of scale invariant maps. The first study of time consistency of dynamic performance measures is due to  Bielecki et al. \cite{BCZ2010}, where the authors elevated the theory of coherent acceptability indices to a dynamic setup in discrete time. It was pointed out  that none of the forms of time consistency for risk measures is suitable for the acceptability indices.

One of the specific features of the acceptability indices, that needed to be accounted for in study of their time consistency, is that these measures of performance can take infinite value. In particular, this required extending the analysis beyond the $L^{\infty}$ framework.

Consequently,  one of the main challenges was to find an appropriate form of time consistency of acceptability indices, that would be  both financially reasonable and mathematically tractable. For the case of random variables (terminal cash flows), the proposed form of time consistency for a dynamic coherent acceptability index $\alpha=\{\alpha_t\}_{0\leq t\leq T}$ reads as follows: for any $\cF_t$-measurable random variables $m_t, \ n_t$, and any $t<T$, the following implications hold
\begin{align}
 \alpha_{t+1}(X)\geq m_t &\quad \Rightarrow \quad \alpha_t(X)\geq m_t, \nonumber \\
 \alpha_{t+1}(X)\leq n_t &\quad \Rightarrow \quad \alpha_t(X)\leq n_t. \label{eq:timeConstAlpha1}
\end{align}
The financial interpretation is clear -- if tomorrow $X$ is acceptable at least at level $m_t$, then today $X$ is also acceptable at least at level $m_t$; similar interpretation holds true for the second part of \eqref{eq:timeConstAlpha1}. It is fair to say, we think, that dynamic acceptability indices and their time consistency properties play a critical role in so called conic approach to valuation and hedging of financial contracts; see  Bielecki et al.\cite{BCIR2012} and and Rosazza Gianin and Sgarra \cite{RosazzaGianinSgarra2012}.

We recall that both risk measures and performance measures, in the nutshell, put preferences on the set of cash flows. While the corresponding forms of time consistency \eqref{eq:strong1} and \eqref{eq:timeConstAlpha1}  for these classes of maps, as argued above, are different, we note that generally speaking both forms of time consistency are linking preferences between different times. The aim of this paper is to present a unified and flexible framework for time consistency of risk measures and performance measures, that integrates existing forms of time consistency.

We consider a (large) class of maps that are postulated to satisfy only two properties - monotonicity and locality\footnote{See Section~\ref{sec:prelims} for rigorous definitions along with a detailed discussion of each property.} - and we study time consistency of such maps. {We focus on these two properties, as, in our opinion, they  have to be satisfied by any reasonable dynamic risk measure or dynamic performance measure.} We introduce the notion of an update rule that is meant to link preferences between different times. The time consistency is defined in terms of an update rule. {It needs to be stressed  that our notion of the update rule is different from the notion of update rule used in  Tutsch \cite{Tutsch2008}. It should be also noted that there exist a large literature in economics, which is focused on the evolution of preferences and the term {\it update rule} is used there as well; see e.g.  Epstein and Schneider \cite{EpsteinSchneider2003}, Hanany and Klibanoff~\cite{HananyKlibanoff2009}, and the references therein. We want to underline that while the concepts of locality and monotonicity are considered there, the updating is applied directly to preference relations, rather than to risk measures or acceptability indices. So, their study is of different nature and is not directly connected to our (axiomatic) framework.}

This paper is the first step in our research leading  towards a unified theory of time consistency of dynamic risk/performance measures {and it ought to be seen as the theoretical basis}. Accordingly, here we focus on formulating and studying the methodological framework without engaging into in-depth presentation of broader aspects of our theory. We refer the reader to our {survey} paper  Bielecki et al. \cite{BCP2014}, where we provide a comprehensive literature overview, present various examples of dynamic LM-measures, update rules and  different types of time consistency, such as middle time consistency, strong time consistency, supermartingale time consistency etc. Moreover, in the survey paper we use our methodology to study connections between these different types of time consistency.

Nevertheless, we spent some time in this paper on illustration of the applicability of our approach. Specifically, we show that almost all known concepts of weak time consistency can be reproduced and studied in terms of a single concept of an update rule introduced in this paper, which is suitable both for dynamic risk measures and dynamic performance measures. In particular, in Proposition~\ref{pr:weak} we characterize weak time consistency for random variables and in Proposition~\ref{pr:proc.weak.semi} we provide a characterization of (semi-)weak time consistency for stochastic processes. Moreover, Propositions~\ref{prop:DCRMtoDAI} and~\ref{prop:DAItoDCRM} show how the weak time consistency property transfers between dynamic coherent risk measures and (normalized) dynamic acceptability indices. This generalizes the result from  Bielecki et al. \cite{BCZ2010} and complements the characterizations from Cherny and Madan~\cite{ChernyMadan2009}, showing how the duality theorems look like in the dynamic setting.

We believe that the general approach introduced in this paper unifies and simplifies the study of time consistency.  A good example of this is Proposition~\ref{prop:mart.ref} that provides a characterization of time consistency via a version of the dynamic programming principle. While in our framework such characterization  is almost immediate, it is not that straightforward to derive it using the benchmark set approach introduced in Tutsch \cite{Tutsch2008}.  Another good example is Proposition~\ref{def:TC-benchSet} where we show how to recover all known (benchmark set) concepts of time consistency using appropriate update rules.

Finally, we want to mention that, traditionally, the investigation of dynamic risk measures and dynamic performances indices is accompanied by robust representation type results. This aspect of the theory is beyond the scope of this study given the generality of the classes of measures considered here.  In particular, the reason for absence in the paper of results regarding robust representation is that such results are  usually derived in the context of convex analysis by exploring convexity (of risk measures) or quasi-concavity (of acceptability indices) properties of some relevant functions. However, we study time consistency without using convex analysis, and we consider functions that are only local and monotone.

The importance of the contribution of the paper can be summarized as follows:
\begin{itemize}
\item We provide a theoretical framework for analysis and synthesis of the various forms of time consistency, allowing for a comparative study of them. Such study is done in the companion paper \cite{BCP2014}.
\item Our theoretical framework is based on the appropriate concept of an update rule. Although the term ``update rule'' has been used in the literature before, the concept of an update rule introduced here is novel and specifically suited for our needs.
\item Our theoretical framework requires minimal assumptions: locality and monotonicity of the measures, for which time consistency is defined and studied.
\end{itemize}

The paper is organized as follows. In Section~\ref{sec:prelims} we introduce some necessary notations and present the main object of our study -- the Dynamic LM-measure.  In Section~\ref{sec:TimeCons} we set forth the main concepts of the paper -- the notion of an updated rule and the definition of time consistency of a dynamic LM-measure. We prove a general result about time consistency, that can be viewed as counterpart of dynamic programming principle \eqref{eq:DPP1}. Additionally, we show that there is a close relationship between update rule approach to time consistency and the approach based on so called benchmark sets. Section~\ref{S:types} is devoted to weak time consistency. The theory presented herein hinges on some new technical results about conditional expectation and conditional essential infimum/supremum for random variables that may take the values $\pm\infty$. These results are presented in Appendix~\ref{A:cond}. To ease the exposition of the main concepts, all technical proofs are deferred to the Appendix~\ref{A:proofs}, unless stated otherwise directly below the theorem or proposition.

\section{Preliminaries}\label{sec:prelims}

Let $(\Omega,\mathcal{F},\bF=\{\mathcal{F}_{t}\}_{t\in\mathbb{T}} ,P)$ be a filtered probability space, with $\mathcal{F}_{0}=\{\Omega,\emptyset\}$,  and $\bT=\set{0,1,\ldots, T}$, for a fixed and finite time horizon $T\in\bN$.\footnote{Most of the results hold true or can be adjusted respectively, to the case of infinite time horizon. For sake of brevity, we will omit the discussion of this case here. }

For $\cG\subseteq\cF$ we denote by $L^0(\Omega,\cG,P)$, and $\bar{L}^0(\Omega,\cG,P)$ the sets of all $\cG$-measurable random variables with values in $(-\infty,\infty)$, and $[-\infty,\infty]$, respectively.
In addition, we will use the notation $L^{p}(\cG):=L^{p}(\Omega,\cG,P)$, $L^{p}_{t}:=L^{p}(\mathcal{F}_{t})$, and $L^{p}:=L_T^{p}$, for $p\in \set{0,1,\infty}$.
Analogous definitions will apply to $\bar{L}^{0}$.
 We will also use the notation $\bV^{p}:=\{(V_{t})_{t\in\bT}: V_{t}\in L^{p}_{t}\}$, for $p\in \{0,1,\infty\}$.

Throughout this paper, $\cX$ will denote either the space of random variables $L^{p}$, or the space of adapted processes $\bV^p$, for $p\in\set{0,1,\infty}$.
If $\cX=L^{p}$,  for $p\in\{0,1,\infty\}$, then the  elements $X \in\cX$ are interpreted as discounted terminal cash flows. On the other hand, if $\cX=\bV^{p}$, for $p\in\{0,1,\infty\}$, then the elements of $\cX$, are interpreted as discounted dividend processes.
It needs to be remarked, that all concepts developed for $\cX=\bV^{p}$ can be easily adapted to the case of cumulative discounted value processes.
The case of random variables can be viewed as a particular case of stochastic processes by considering cash flows with only the terminal payoff, i.e. stochastic processes such that $V=(0,\ldots,0,V_T)$. Nevertheless, we treat this case separately for transparency.
For both cases  we will consider standard pointwise order, understood in the almost sure sense.
In what follows, we will also make use of the multiplication operator denoted as $\cdot_{t}$ and defined by:
\begin{align}
m\cdot_{t}V &:=(V_{0},\ldots,V_{t-1},mV_{t},mV_{t+1},\ldots), \nonumber\\
m\cdot_{t}X &:= mX,\label{eq:conventionV}
\end{align}
for  $V\in\Set{(V_t)_{t\in\bT} \mid V_{t}\in L^{0}_{t}}$, $X\in L^{0}$ and  $m\in L^{\infty}_{t}$.
In order to ease the notations, if no confusion arises, we will drop $\cdot_t$ from the above product, and we will simply write $mV$ and $mX$ instead of $m\cdot_{t}V$ and $m\cdot_{t}X$, respectively.

\begin{remark}
We note that the space $\bV^{p}, \ p\in\{0,1,\infty\}$, endowed with  multiplication $(\,\cdot_{t},)$ does not define a proper $L^{0}$--module  \cite{FilipovicKupperVogelpoth2009} (e.g. $0\cdot_{t} V\neq 0$ for some $V\in\bV^{p}$). However, in what follows, we will adopt some concepts from $L^0$-module theory which naturally fit into our study.
Moreover, as in many cases we consider, if one additionally assumes {\it independence of the past}, and replaces $V_{0},\ldots,V_{t-1}$ with $0$s in \eqref{eq:conventionV}, then $\cX$ becomes an $L^{0}$--module. We refer the reader to \cite{BCDK2013,BCP2013} for a thorough discussion on this matter.
\end{remark}

Throughout, we will use the convention that $\infty-\infty=-\infty+\infty=-\infty$ and $0\cdot\pm\infty=0$.

\noindent
For $t\in\bT$ and $X\in\bar{L}^{0}$, we define the (generalized) $\cF_{t}$-conditional expectation of $X$ by
$$
E[X|\cF_{t}]:=\lim_{n\to\infty}E[(X^{+}\wedge n)|\cF_{t}]-\lim_{n\to\infty}E[(X^{-}\wedge n)|\cF_{t}],
$$
where $X^{+}=(X\vee 0)$ and $X^{-}=(-X\vee 0)$. {Note  that, in view of our convention we have that $(-1)(\infty-\infty) = \infty \neq -\infty + \infty =-\infty$, which, in particular, implies that we might get $-E[X]\neq E[-X]$.}
{Thus, the conditional expectation operator defined above is no longer linear on $\bar{L}^{0}$ space (see Appendix~\ref{A:cond}, Proposition \ref{pr:condexp}). Similarly, for any $t\in\bT$ and $X\in\bar{L}^{0}$, we define the (generalized) $\cF_{t}$-conditional essential infimum by\footnote{Since both sequences $\Essinf_{t}(X^{+}\wedge n)$ and $\Esssup_{t}(X^{-}\wedge n)$ are monotone, the corresponding limits exist.}
\begin{equation}\label{eq:essinf}
\Essinf_{t}X:=\lim_{n\to\infty}\Big[\Essinf_{t}(X^{+}\wedge n)\Big]-\lim_{n\to\infty}\Big[\Esssup_{t}(X^{-}\wedge n)\Big],
\end{equation}
and respectively, we put $\Esssup_{t}(X):=-\Essinf_{t}(-X)$. For some basic properties of this operator and the definition of conditional essential infimum on $L^{\infty}$ see Appendix~\ref{A:cond}. In particular, note that, for any $X\in \bar{L}^{0}_{t}$, we get $\Essinf_{t}X=X$.

Next, we introduce the main object of this study.
\begin{definition}
A family $\varphi=\{\varphi_{t}\}_{t\in\mathbb{T}}$ of maps $\varphi_{t}:\cX\to\bar{L}^{0}_{t}$ is a {\it Dynamic LM-measure} if $\varphi$ satisfies
\begin{enumerate}[1)]
\item {\it (Locality)}  $\1_{A}\varphi_{t}(X)=\1_{A}\varphi_{t}(\1_{A}\cdot_{t} X)$;
\item {\it (Monotonicity)} $X\leq Y \Rightarrow \varphi_{t}(X)\leq \varphi_{t}(Y)$;
\end{enumerate}
for any $t\in\bT$, $X,Y\in\cX$, and $A\in\cF_{t}$.
\end{definition}

We believe that locality and monotonicity are two properties that must be satisfied by any reasonable dynamic measure of performance and/or measure of risk. Monotonicity property is natural for any numerical representation of an order between elements of $\cX$.  {The locality property essentially means that the values of the LM-measure restricted to a set $A\in\cF$ remain invariant with respect to the values of the arguments outside of the same set $A\in\cF$; in particular, the events that will not happen in the future do not change the value of the measure today.  }

\label{page:footnote} {Dynamic LM-measures contain several important subclasses.} Among the most recognized ones are {dynamic risk measures} and dynamic performance measures {(dynamic acceptability indices)}. These classes of measures have been extensively studied in the literature over the past decade.

We recall that a function $\varphi_{t}:\cX\to\bar{L}^{0}_{t}$ is: \textit{cash additive} if $\varphi(X+m1_{\{t\}})=\varphi_t(X)+m$, for any $X\in\cX$, $t\in\bT$, and $m\in L_t^{p}$; \textit{scale invariant} if $\varphi_t(\beta\cdot_{t}X)=\varphi_t(X)$, for any $X\in\cX$, $t\in\bT$, and $\beta\in L_t^p, \beta >0$.

\textit{A dynamic monetary utility measure} is a cash-additive LM-measure, and a dynamic risk measure is the negative of a dynamic monetary utility measure. For convenience, we will study dynamic monetary utility measure in this study rather than dynamic risk measures. Cash additivity is the key property that distinguishes utility/risk measures from all other measures. This property means that adding $\$m$ to a portfolio today reduces the overall risk by the same amount $\$m$. From the regulatory perspective, the value of a risk measure is typically interpreted as the minimal capital requirement for a bank. For more details on coherent/covex/monetary risk measures we refer the reader to the survey papers \cite{FollmerSchied2010,AcciaioPenner2010}.

\textit{A dynamic performance measure} is a scale invariant LM-measure. As already mentioned, the distinctive property of performance measures is the scale invariance - a rescaled portfolio or a cash flow is accepted at the same level. Performance measures, sometimes referred to as acceptability indices,  were studied in~\cite{ChernyMadan2009,BCZ2010,CheriditoKromer2012,BCP2013}, and they are meant to provide an assessment of how good a financial position is.\footnote{Some authors treat acceptability indices as the special subclass of performance measures, that satisfy the quasi-concavity axiom. In particular,~\cite{CheriditoKromer2012} gives examples of performance indices that are not quasi-concave. Nevertheless, in this paper we have decided to use those two names interchangeably.}
It needs to be noted that the theory developed in this paper can also be applied to sub-scale invariant dynamic assessment indices  studied in \cite{RosazzaGianinSgarra2012,BCC2014}.

\section{Time consistency and update rules}\label{sec:TimeCons}
In this section we introduce the main concept of this paper - the time consistency of {dynamic risk measures and dynamic performance measures}, or more generally, the time consistency of dynamic LM-measures introduced in the previous section.

We recall that these dynamic LM-measures are defined on $\cX$, where  $\cX$ either denotes the space  $L^{p}$ of random variables or the space $\bV^{p}$ of stochastic processes, for $p\in\set{0,1,\infty}$, so, our study of time consistency is done relative to such spaces. Nevertheless, the definition of time consistency can be easily adapted to more general spaces, such as Orlicz hearts (as studied in \cite{CheriditoLi2009}), or, such as topological $L^{0}$-modules (see for instance~\cite{BCDK2013}).

Assume that   $\varphi$ is a dynamic LM-measure on $\cX$.
For  an arbitrary fixed $X\in\cX$ and $t\in \bT$, the value $\varphi_{t}(X)$ represents a quantification (measurement) of preferences about $X$ at time $t$. Clearly, it is reasonable to require that any such quantification (measurement) methodology should be coherent as time passes. This is precisely the motivation behind the concepts of time consistency of dynamic LM-measures.

{There are various forms of time consistency proposed in the literature, some of them suitable for one class of measures, others for a different class of measures.}
For example, for dynamic convex (or coherent) risk measures, various version of time consistency surveyed in \cite{AcciaioPenner2010} can be {seen} as versions of the celebrated  dynamic programming principle.
On the other hand, as shown in \cite{BCZ2010}, dynamic programming principle essentially is not suited for scale invariant measures such as dynamic acceptability indices, and the authors introduce a new type of time consistency, tailored for these measures, and provide a robust representation of them.
Nevertheless, in all these cases the time consistency property connects, in a {noncontradictory}  way, the measurements {done} at different times.

Next, we will introduce the notion of update rule that serves as the main tool in relating the measurements of preferences at different times, and also, it is the main building block of our unified theory of time consistency property.

\begin{definition}\label{def:UMST}
We call a family $\mu=\{\mu_{t,s}:\, t,s\in\bT,\, s>t\}$ of maps $\mu_{t,s}:\bar{L}^{0}_{s}\times\cX\to\bar{L}^{0}_{t}$ an {\it update rule} if for any $s>t$, the map $\mu_{t,s}$ satisfies the following properties:
\begin{enumerate}[1)]
\item (Locality) $\1_{A}\mu_{t,s}(m,X)=\1_{A}\mu_{t,s}(\1_{A}m,X)$;
\item (Monotonicity) if $m\geq m'$, then  $\mu_{t,s}(m,X)\geq \mu_{t,s}(m',X)$;
\end{enumerate}
for any $X\in\cX$, $A\in\cF_{t}$ and $m,m'\in\bar{L}^{0}_{s}$.
\end{definition}
Since LM-measures are local and monotone, properties with clear financial interpretations, the update rules are naturally assumed to be local and monotone too.

The first argument $m\in\bar{L}^{0}_{s}$ in $\mu_{t,s}$  serves as a benchmark to which the measurement $\varphi_{s}(X)$ is compared. The presence of the second argument, $X\in\cX$, in $\mu_{t,s}$, allows the update rule to depend on the objects to which the preferences are applied to. However, as we will see in next section, there are natural situations when the update rules are independent of $X\in\cX$, and sometimes they do not even depend on the future times $s\in\bT$.

{
\begin{remark}\label{rem:super1}
As we have mentioned, the update rule is used for updating preferences through time. This, for example, can be achieved in terms of the conditional expectation operator
\begin{equation}\label{eq:ur:exp}
\mu_{t,s}(m,X)=E[m|\cF_{t}],
\end{equation}
which is an update rule.
Note that this particular update rule does not depend on $s$ or $X$. Update rules might be also used for discounting the preferences. Intuitively speaking, the risk of loss in the far future might be more preferred than the imminent risk of loss (see~\cite{Cherny2010} for the more detailed explanation of this idea). For example, the update rule $\mu$ of the form
\begin{equation}\label{eq:ur:exp3}
\mu_{t,s}(m,X)=\left\{
\begin{array}{ll}
\varepsilon^{s-t} E[m|\cF_{t}] & \textrm{on } \{E[m|\cF_{t}] \geq 0\},\\
\varepsilon^{t-s} E[m|\cF_{t}] &  \textrm{on } \{E[m|\cF_{t}] < 0\}.
\end{array}\right.
\end{equation}
for a fixed $\varepsilon\in (0,1)$ would achieve this goal. Note that `discounting' proposed here has nothing to do with the ordinary discounting, as we act on discounted values already.
\end{remark}
}
Next, we define several particular classes of update rules, suited for our needs.
\begin{definition}\label{UM}
Let $\mu$ be an update rule. We say that $\mu$ is:
\begin{enumerate}[1)]
\item {{\it $X$-invariant}, if $\mu_{t,s}(m,X)=\mu_{t,s}(m,0)$;}
 \item {\it $sX$-invariant}, if there exists a family $\{\mu_{t}\}_{t\in\bT}$ of maps $\mu_{t}:\bar{L}^{0}\to\bar{L}^{0}_{t}$, such that $\mu_{t,s}(m,X)=\mu_{t}(m)$;
 \item {\it Projective}, if it is $sX$-invariant and $\mu_{t}(m_{t})=m_{t}$;
 \end{enumerate}
for any $s,t\in\bT$, $s>t$, $X\in\cX$, $m\in\bar{L}^{0}_{s}$ and $m_{t}\in\bar{L}^{0}_{t}$.
\end{definition}
Examples of update rules satisfying 1) and 3) are given by \eqref{eq:ur:exp3} and \eqref{eq:ur:exp}, respectively. The update rule, which satisfies 2), but not 3) can be constructed by substituting $\varepsilon^{t-s}$ with a constant in \eqref{eq:ur:exp3}. Generally speaking update rules for stochastic processes will not satisfy 1) as the information about the process on the time interval $(t,s)$ will affect $\mu_{t,s}$; see Subsection~\ref{S:semi.weak} for details.

\begin{remark}\label{rem:UMnotation}
If an update rule $\mu$ is $sX$-invariant, then  it is enough to consider only the corresponding family $\{\mu_{t}\}_{t\in\bT}$. Hence, with slight abuse of notation we will write $\mu=\{\mu_{t}\}_{t\in\bT}$, and call it an update rule as well.
\end{remark}

We are now ready to introduce the general definition of time consistency.

\begin{definition}$\!\!\!$\footnote{We introduce the concept of time consistency only for LM-measures, as this is the only class of measures used in this paper. However, the definition itself is suitable for any map acting from $\cX$ to $\bar{L}^0$. For example, traditionally in the literature, the time consistency is defined for dynamic risk measures (negatives of cash-additive LM-measures), and the above definition of time consistency will be appropriate, although one has to flip `acceptance' with `rejection'.
}
Let $\mu$ be an update rule. We say that the dynamic LM-measure $\varphi$ is {\it $\mu$-acceptance (resp. $\mu$-rejection) time consistent} if
\begin{equation}\label{eq:atc}
\varphi_{s}(X)\geq m_{s}\quad (\textrm{resp.} \leq) \quad \Longrightarrow\quad \varphi_{t}(X)\geq \mu_{t,s}(m_{s},X)\quad (\textrm{resp.} \leq),
\end{equation}
for all $s,t\in\bT$, $s>t$, $X\in \cX$ and $m_{s}\in \bar{L}^{0}_{s}$.
If property \eqref{eq:atc} is satisfied only for $s,t\in\bT$, such that $s=t+1$, then we say that $\varphi$ is {\it one step $\mu$-acceptance (resp. one step $\mu$-rejection) time consistent}.
\end{definition}

The financial interpretation of acceptance time consistency is straightforward: if $X\in\cX$ is accepted at some future time $s\in\bT$, at least at level $m$, then today, at time $t\in\bT$, it is accepted at least at level $\mu_{t,s}(m,X)$. Similarly for rejection time consistency. Essentially, the update rule $\mu$ translates the preference levels at time $s$ to preference levels at time $t$. As it turns out, this simple and intuitive definition of time consistency, with appropriately chosen $\mu$, will cover various cases of time consistency for risk and performance measures that can be found in the existing literature (see \cite{BCP2014} for a survey).

Next, we will give an equivalent formulation of time consistency, which, in fact, might be taken as a definition of time consistency (in place of \eqref{eq:atc}). Given the nature of the update rule and its purpose, we however believe that property \eqref{eq:atc} is more natural defining property, as compared to \eqref{eq:accepTimeConsAlt}. While the proof of the equivalence is simple, the result itself is very important and it will be conveniently  used in the sequel. Moreover, it can be viewed as a counterpart of dynamic programming principle, which is an equivalent formulation of dynamic consistency for convex/coherent risk measures. This is the reason why we separate out this result in the form of proposition.

\begin{proposition}\label{prop:mart.ref}
Let $\mu$ be an update rule, and let $\varphi$ be a dynamic LM-measure. Then, $\varphi$ is $\mu$-acceptance  (resp. $\mu$-rejection) time consistent if and only if
\begin{equation}\label{eq:accepTimeConsAlt}
\varphi_{t}(X)\geq \mu_{t,s}(\varphi_{s}(X),X)\quad (\textrm{resp.} \leq),
\end{equation}
for any $X\in\cX$ and $s,t\in\bT$, such that $s>t$.
\end{proposition}

\begin{remark}\label{rem:preservetc}
It is clear, and also naturally desired, that a monotone transformation of an LM-measure will not change the preference order of the underlying elements.
We want to emphasize  that a monotone transformation will also preserve the time consistency. In other words, the preference orders will be also preserved in time.
Indeed, if $\varphi$ is $\mu$-acceptance time consistent, and $g:\bar{\bR} \to\bar{\bR}$ is a strictly monotone function, then the family $\{g\circ \varphi_{t}\}_{t\in\mathbb{T}}$ is $\widetilde{\mu}$-acceptance time consistent, where the update rule $\widetilde{\mu}$ is defined by $\widetilde{\mu}_{t,s}(m,X)=g(\mu_{t,s}(g^{-1}(m),X))$, for $t,s\in\bT$, $s>t$, $X\in\cX$ and $m\in \bar{L}^{0}_{s}$.
\end{remark}

In the case of random variables, $\cX=L^p$,  we we will usually consider update rules that are $X$-invariant. The case of stochastic processes is more intricate. If $\varphi$ is a dynamic LM-measure, and $V\in\bV^p$, then in order to compare $\varphi_{t}(V)$ and $\varphi_{s}(V)$, for $s>t$, one also needs to take into account the cash flows between times $t$ and $s$. Usually, for $\cX=\bV^{p}$ we consider update rules, such that
\begin{equation}\label{eq:prTOrv.f}
\mu_{t,t+1}(m,V)=\mu_{t,t+1}(m,0)+f(V_{t}),
\end{equation}
where $f:\bar{\bR}\to \bar{\bR}$ is a Borel measurable function, such that $f(0)=0$. We note, that any such one step update rule $\mu$ can be easily adapted to the case of random variables. Indeed, upon  setting $\widetilde{\mu}_{t,t+1}(m):=\mu_{t,t+1}(m,0)$ we get a one step $X$-invariant update rule $\widetilde{\mu}$, which is suitable for random variables. Moreover,   $\widetilde{\mu}$ will define the corresponding type of one step time consistency for random variables. Of course, this correspondence between update rule for processes and random variables is valid only for `one step' setup.

Moreover, for update rules, which admit the so called nested composition property~(cf. \cite{Ruszczynski2010,RuszczynskiShapiro2006a} and references therein),
\begin{equation}\label{eq:RvToPr}
{\mu_{t,s}(m,V)=\mu_{t,t+1}(\mu_{t+1,t+2}(\ldots\mu_{s-2,s-1}(\mu_{s-1,s}(m,V),V)\ldots V),V),}
\end{equation}
we have that $\mu$-acceptance (resp. $\mu$-rejection) time consistency is equivalent to one step $\mu$-acceptance (resp. $\mu$-rejection) time consistency.

\subsection{Relation between update rule approach and the benchmark  approach}
As we will show in this section, there is a close relationship between our update rule approach to time consistency and the approach based on so called benchmark sets. The latter approach was initiated in \cite{Tutsch2008}, where the author applied it in the context of dynamic risk measures. Essentially, a  benchmark set is a collection of elements from $\cX$ that satisfy some additional structural properties.

For simplicity, we shall assume here that $\cX=L^p$, for $p\in\{0,1,\infty\}$. The definition of time consistency in terms of benchmark  sets is as follows:

\begin{definition}\label{def:TC-benchSet}
Let $\varphi$ be a dynamic LM-measure and let $\cY=\{\cY_{t}\}_{t\in\bT}$ be a {\it family of benchmark sets}, that is, sets $\cY_t$ such that $\cY_{t}\subseteq L^p$, $0\in\cY_{t}$ and $\cY_t+\bR=\cY_t$. We say that $\varphi$ is {\it acceptance (resp. rejection) time consistent with respect to $\cY$}, if
\begin{equation}\label{eq:benchmark2}
\varphi_{s}(X)\geq \varphi_{s}(Y)\quad (resp. \leq)\quad \Longrightarrow\quad \varphi_{t}(X)\geq \varphi_{t}(Y)\quad (resp. \leq),
\end{equation}
for all $s\geq t$, $X\in L^p$ and $Y\in\cY_{s}$.
\end{definition}
Informally, the `degree' of time consistency with respect to $\cY$ is measured by the size of $\cY$.
Thus, the larger the sets $\cY_{s}$ are, for each $s\in\bT$, the stronger is the degree of time consistency of $\varphi$.

We now have the following important proposition,

\begin{proposition}\label{prop:benchmark}
Let $\varphi$ be a dynamic LM-measure and let $\cY$ be a family of benchmark sets.
 Then, there exists an update rule $\mu$ such that: $\varphi$ is acceptance (resp. rejection) time consistent with respect to $\cY$ if and only if it is $\mu$-acceptance (resp. $\mu$-rejection) time consistent.
\end{proposition}

\noindent The update rule $\mu$ is said to provide $\varphi$ with the same type of time consistency as $\cY$ does, and vice versa. Generally speaking, the converse implication does not hold true, i.e. given an LM-measure $\varphi$ and an update rule $\mu$ it may not be possible to construct $\cY$ so that it provides the same type of time consistency as $\mu$ does. In other words, the notion of time consistency given in terms of updates rule is more general.

\section{Weak time consistency}\label{S:types}
In this section we will discuss examples of update rules, which relate to weak time consistency for random variables and for stochastic processes. This is meant to illustrate the framework developed earlier in this paper. As mentioned in the Introduction,  see~\cite{BCP2014} for a comprehensive survey of various types of time consistency and connections between them.


The notion of weak time consistency was introduced in~\cite{Tutsch2008}, and subsequently studied in \cite{AcciaioPenner2010,ArtznerDelbaenEberHeathKu2007,CheriditoDelbaenKupper2006,DetlefsenScandolo2005,
AcciaioFollmerPenner2010}.
The idea is that if `tomorrow', say at time $s$, we accept $X\in\cX$ at level $m_s\in\cF_s$, then `today', say at time  $t$, we would accept  $X$ at least at any level {lower} or equal {to} $m_s$, {appropriately} adjusted by the information $\cF_t$ available at time $t$ (cf. \eqref{eq:inft}). Similarly, if tomorrow we reject $X$ at level {higher or equal to} $m_s\in\cF_s$, then today, we should also reject $X$ at any level {higher} than $m_s$, adjusted to the flow of information $\cF_t$. This suggests that the update rules should be taken as  $\cF_t$-conditional essential infimum and supremum, respectively. Towards this end, we first show that $\cF_t$-conditional essential infimum and supremum are projective update rules.

\begin{proposition}\label{th:essinf}
The family $\mu^{\inf}:=\{\mu_{t}^{\inf}\}_{t\in\bT}$ of maps $\mu^{\inf}_{t}:\bar{L}^{0}\to\bar{L}^{0}_{t}$ given by
\[
\mu^{\inf}_{t}(m)=\Essinf_{t}m,
\]
is a projective\footnote{See Remark~\ref{rem:UMnotation} for the comment about notation.\label{fn:notation}} update rule. Similar result is true for the family $\mu^{\sup}:=\{\mu_{t}^{\sup}\}_{t\in\bT}$ of maps $\mu^{\sup}_{t}:\bar{L}^{0}\to\bar{L}^{0}_{t}$ given by $\mu^{\sup}_{t}(m)=\Esssup_{t}m$.
\end{proposition}

\subsection{Weak time consistency for random variables}
Recall that the case of random variables corresponds to $\cX=L^{p}$, for a fixed $p\in\{0,1,\infty\}$.
We proceed with the definition of weak acceptance and weak rejection time consistency (for random variables).

\begin{definition}\label{type.of.cons.weak}
A dynamic LM-measure $\varphi$ is said to be {\it weakly acceptance (resp. {\it  weakly rejection}) time consistent}  if it is $\mu^{\inf}$-acceptance (resp. $\mu^{\sup}$-rejection) time consistent.
\end{definition}

Definition~\ref{type.of.cons.weak}  of time consistency is equivalent to many forms of time consistency studied in the current literature.
Usually,  the weak time consistency is considered for dynamic monetary risk measures on $L^{\infty}$ (cf.~\cite{AcciaioPenner2010} and references therein); we refer to this case  as to} the `classical weak time consistency.'
It was {observed }in \cite{AcciaioPenner2010} that in the classical weak time consistency framework, weak acceptance (resp. weak rejection) time consistency  is equivalent to the statement that for any $X\in L^p$ and $s>t$, we get
\begin{equation}\label{eq:classicWeak}
\varphi_{s}(X)\geq 0 \Rightarrow \varphi_{t}(X)\geq 0,\quad\quad \textrm{(resp. $\leq$)}.
\end{equation}
This observation was the motivation  for our definition of weak acceptance (resp. weak rejection) time consistency, and the next proposition explains why so.

\begin{proposition}\label{pr:weak}
Let $\varphi$ be a dynamic LM-measure. The following conditions are equivalent
\begin{enumerate}[1)]
\item $\varphi$ is weakly acceptance time consistent.
\item For any $X\in  L^p$, $s,t\in\bT$, $s>t$, and $m_{t}\in \bar{L}^{0}_{t}$,
\[
\varphi_{s}(X)\geq m_{t} \Rightarrow \varphi_{t}(X)\geq m_{t}.
\]
\end{enumerate}
If additionally $\varphi$ is a normalized dynamic monetary utility measure\footnote{i.e $\varphi_{t}(0)=0$ and $\varphi_{t}(X+\
c_{t})=\varphi_{t}(X)+c_{t}$ for any $t\in\bT$, $X\in  L^p$ and $c_{t}\in L^{\infty}_{t}$.}, then the above conditions are equivalent to
\begin{enumerate}[1)]
\item[3)] For any $X\in  L^p$ and $s,t\in\bT$, $s>t$,
\[
\varphi_{s}(X)\geq 0 \Rightarrow \varphi_{t}(X)\geq 0.
\]
\end{enumerate}
Similar result holds true for the weak rejection time consistency.
 \end{proposition}

Property 2) in Proposition~\ref{pr:weak} was also suggested as the notion of (weak) acceptance and (weak) rejection time consistency in the context of scale invariant measures (cf.~\cite{BiaginiBion-Nadal2012,BCZ2010}).

In many papers studying risk measurement theory (cf. \cite{DetlefsenScandolo2005} and references therein), the weak form of time consistency is defined using dual approach to the measurement of risk. Rather than directly updating the level of preferences $m$, as in our approach, in the dual approach the level of preference is updated indirectly by manipulating probabilistic scenarios and explaining the update procedure by using so called {\it pasting property}  (see e.g.~\cite[Def. 9]{DetlefsenScandolo2005}). As shown in the next result, our update rule related to weak form of time consistency admits dual representation, allowing us to link our definition with the dual approach.

\begin{proposition}\label{pr:essinf.cond.aa}
For any $m\in\bar{L}^{0}$ and $t\in\bT$, we have
\begin{equation}\label{eq:essinf.eq.pt}
\mu^{\inf}_{t}(m)=\essinf_{Z\in P_{t}}E[Zm|\mathcal{F}_{t}],
\end{equation}
where $P_{t}:=\{Z\in L^{0} \mid Z\geq 0,\ E[Z|\mathcal{F}_{t}]=1\}$. Similar result is true for $\Esssup_{t}m$.
\end{proposition}
In \eqref{eq:essinf.eq.pt}, the random variables $Z\in P_{t}$ may be treated as the Radon-Nikodym derivatives with resect to $P$  of some probability measures $Q$ such that $Q\ll P$ and $Q|_{\cF_{t}}=P|_{\cF_{t}}$. {The family $P_{t}$ may thus be thought of as the family of all possible $\cF_{t}$-conditional probabilistic scenarios. Accordingly, $\mu^{\inf}_{t}(m)$ represents the $\cF_{t}$-conditional worst-case preference update with respect to all such scenarios.} Note that by combining Propositions~\ref{prop:mart.ref} and \ref{pr:essinf.cond.aa}, we obtain that weak acceptance time consistency of $\varphi$ is equivalent to the condition
\begin{equation}\label{eq:robust.expl}
\varphi_{t}(X) \geq \essinf_{Z\in P_{t}}E[Z\varphi_{s}(X)|\mathcal{F}_{t}],
\end{equation}
which  is a starting point for almost all robust definitions of weak time consistency, for $\varphi$'s that admit dual representation~\cite{DetlefsenScandolo2005}.

As next result shows, the weak time consistency is indeed one of the weakest forms of time consistency, being implied by any other concept of  time consistency generated by a projective rule.

\begin{proposition}\label{pr:UDMprop}
Let $\varphi$  be a dynamic LM-measure and let $\mu$ be a projective update rule. If $\varphi$ is $\mu$-acceptance (resp. $\mu$-rejection) time consistent, then $\varphi$ is weakly acceptance (resp. weakly rejection) time consistent.
\end{proposition}

In particular, recall that time consistency is preserved under monotone transformations, Remark~\ref{rem:preservetc}. Thus, for any strictly monotone function $g:\bar{\bR}\to\bar{\bR}$ , if $\varphi$ is weakly acceptance (resp. weakly rejection) time consistent, then $\{g\circ\varphi_{t}\}_{t\in\bT}$ also is weakly acceptance (resp. weakly rejection) time consistent.

\subsection{Weak and Semi-weak time consistency for stochastic processes}\label{S:semi.weak}

In this subsection we introduce {and discuss} the concept of semi-weak time consistency for stochastic processes. Thus, we take $\cX=\bV^{p}$, for a fixed $p\in\{0,1,\infty\}$. As it will turn out, in the case of random variables semi-weak time consistency coincides with  weak time consistency; that is why we omitted discussion of semi-weak consistency in the previous section.

To provide a better perspective for the concept of semi-weak time consistency, we start with the definition of weak time consistency for stochastic processes, which transfers directly from the case of random variables by using \eqref{eq:prTOrv.f}.

\begin{definition}\label{type.of.cons.proc.weak}
Let $\varphi$ be a dynamic LM-measure.
We say that $\varphi$ is {\it weakly acceptance (resp. weakly rejection) time consistent for stochastic processes} if it is one step $\mu$-acceptance (resp. one step $\mu^{*}$-rejection) time consistent, where the update rule is given by
$$
\mu_{t,t+1}(m,V)=\mu^{\inf}_{t}(m)+V_{t}\qquad (\textrm{resp. } \mu^{*}_{t,t+1}(m,V)=\mu^{\sup}_{t}(m)+V_{t}).
$$
\end{definition}

As mentioned earlier, the update rule, and consequently weak time consistency for stochastic processes, depends also on the value of the process (the dividend paid)  at time $t$.
If tomorrow, at time $t+1$, we accept $X\in  \bV^p$ at level greater than $m_{t+1}\in\cF_{t+1}$, then today at time $t$, we will accept $X$ at least at level $\Essinf_t m_{t+1}$ (i.e. the worst level of $m_{t+1}$ adapted to the information $\cF_t$) plus the dividend $V_t$ received today.

For counterparts of Propositions~\ref{pr:weak} and~\ref{pr:UDMprop} for the case of stochastic processes see ~\cite{BCP2014}.

As it was shown in \cite{BCZ2010}, none of the existing, at that time, forms of time consistency were suitable for scale invariant maps, such as acceptability indices. In fact, even the weak acceptance and the weak rejection time consistency for stochastic processes are too strong in case of {acceptability indices}.
Because of that, we need  a weaker notion of time consistency, which we will refer to as semi-weak acceptance and semi-weak rejection time consistency. These notions of time consistency are suited precisely for {acceptability indices}, and we refer the reader to \cite{BCZ2010} for a detailed discussion on time consistency for {acceptability indices} and their dual representations\footnote{In \cite{BCZ2010} the authors combined both semi-weak acceptance and rejection time consistency into one single definition and call it time consistency.}.

\begin{definition}\label{def:semi}
Let $\varphi$ be a dynamic LM-measure (for processes). Then $\varphi$ is said to be:
\begin{itemize}
\item {\it Semi-weakly acceptance time consistent}  if it is one step $\mu$-acceptance time consistent, where the update rule is given by
$$
\mu_{t,t+1}(m,V) =1_{\{V_{t}\geq 0\}}\mu^{\inf}_{t}(m)+1_{\{V_{t}< 0\}}(-\infty).\quad
$$
\item {\it Semi-weakly rejection time consistent}  if it is one step $\mu'$-rejection time consistent, where the update rule is given by
$$
\mu'_{t,t+1}(m,V) =1_{\{V_{t}\leq 0\}}\mu^{\sup}_{t}(m)+1_{\{V_{t}> 0\}}(+\infty).
$$
\end{itemize}
\end{definition}

It is straightforward to check that weak acceptance/rejection time consistency for stochastic processes always implies semi-weak acceptance/rejection time consistency.

Next, we will show that the definition of semi-weak time consistency is indeed equivalent to the time consistency introduced in  \cite{BCZ2010}, and later studied in \cite{BiaginiBion-Nadal2012,BCC2014}.

\begin{proposition}\label{pr:proc.weak.semi}
Let $\varphi$ be a dynamic LM-measure on $\bV^p$ . The following conditions are equivalent.
\begin{enumerate}[1)]
\item $\varphi$ is semi-weakly acceptance time consistent, i.e. for all $V\in \bV^p$, $t\in\bT$,  $t<T$, and $m_{t}\in \bar{L}^{0}_{t}$,
$$
\varphi_{t+1}(V)\geq m_{t+1} \Rightarrow \varphi_{t}(V) \geq \1_{\{V_{t}\geq 0\}}\Essinf_t(m_{t+1})+1_{\{V_{t}< 0\}}(-\infty).
$$
\item For all $V\in \bV^p$ and $t\in\bT$,  $t<T$, $\varphi_{t}(V) \geq \1_{\{V_{t}\geq 0\}}\Essinf_t(\varphi_{t+1}(V))+1_{\{V_{t}< 0\}}(-\infty)$.
\item For all $V\in \bV^p$, $t\in\bT$,  $t<T$, and $m_{t}\in \bar{L}^{0}_{t}$, such that $V_{t}\geq 0$ if $\varphi_{t+1}(V)\geq m_{t}$, then $\varphi_{t}(V)\geq m_{t}$.
\end{enumerate}
Similar result is true for semi-weak rejection time consistency.
\end{proposition}
Property 3) in Proposition~\ref{pr:proc.weak.semi} illustrates best the financial meaning of semi-weak acceptance time consistency: if tomorrow we accept the dividend stream $V\in \bV^p$ at level $m_t$,  and if we get a positive dividend $V_t$ paid today at time $t$, then today we accept the cash flow $V$ at least at level $m_t$ as well. Similar interpretation is valid for semi-weak rejection time consistency.

The next two results give an important (dual) connection between cash additive risk measures and {acceptability indices. In particular, these results shed light on the relation between time consistency property of dynamic acceptability indices, represented by the family $\{\alpha_{t}\}_{t\in\bT}$ below,  and time consistency of the corresponding family $\{\phi^{x}\}_{x\in\bR_{+}}$, where $\phi^{x}=\{\phi^{x}_{t}\}_{t\in\bT}$ is a dynamic risk measure (for any $x\in\bR_+$).}

\begin{proposition}\label{prop:DCRMtoDAI}
 Let $\{\phi^{x}\}_{x\in\bR_{+}}$ be a decreasing family of dynamic LM-measures\footnote{A family, indexed by $x\in\bR_{+}$, of maps $\{\phi_{t}^{x}\}_{t\in\bT}$,  will be called {\it decreasing}, if $\phi_{t}^{x}(X)\leq \phi_{t}^{y}(X)$ for all $X\in\cX$, $t\in\bT$ and $x,y\in\bR_{+}$, such that $x\geq y$.}.
    Assume that for each $x\in\bR_{+}$, $\{\phi^{x}_{t}\}_{t\in\bT}$ is weakly acceptance (resp. weakly rejection) time consistent.
    Then, the family $\{\alpha_{t}\}_{t\in\bT}$ of maps $\alpha_{t}:\bV^p\to \bar{L}^{0}_{t}$ defined by\footnote{Note that the map defined in \eqref{eq:DCRMtoDAI} is $\cF_{t}$-measurable as the essential supremum over an  uncountable family of $\cF_{t}$-measurable random variables. See Appendix \ref{A:cond}.}
\begin{equation}\label{eq:DCRMtoDAI}
\alpha_{t}(V):=\esssup_{x\in \bR^{+}}\{x \1_{\{\phi_{t}^{x}(V)\geq0\}}\},
\end{equation}
is a semi-weakly acceptance (resp. semi-weakly rejection) time consistent dynamic LM-measure.
\end{proposition}
Observe that a version of $\alpha_{t}(V)$ is given as
\begin{equation}\label{eq:DCRMtoDAI2}
\alpha_{t}(V)(\omega)=\sup\{x\in\bR_{+}: \phi_{t}^{x}(V)(\omega)\geq0\}.
\end{equation}
As the representation \eqref{eq:DCRMtoDAI2} is more convenient than \eqref{eq:DCRMtoDAI}, it will be used in the proofs given in the Appendix.

\begin{proposition}\label{prop:DAItoDCRM}
Let $\{\alpha_{t}\}_{t\in\bT}$ be a dynamic LM-measure, which is independent of the past and translation invariant\footnote{We say that $\alpha$ is {\it translation invariant}  if $\alpha_{t}(V+m1_{\set{t}})=\alpha_{t}(V+m1_{\set{s}})$, for any $m\in L^{p}_{t}$ and $V\in\bV^p$, where $1_{\set{t}}$ corresponds to process equal to $1$ a time $t$ and 0 elsewhere; We say that $\alpha$ is {\it independent of the past} if $\alpha_{t}(V)=\alpha_t((0,\ldots, 0, V_t,\ldots, V_T))$, for any $V\in\bV^p$.}.
Assume that $\{\alpha_{t}\}_{t\in\bT}$ is semi-weakly acceptance (resp. semi-weakly rejection) time consistent.
Then, for any $x\in\bR_{+}$, the family $\{\phi_{t}^{x}\}_{t\in\bT}$ defined by
\begin{equation}\label{eq:DAItoDCRM}
\phi^{x}_{t}(V)=\essinf_{c\in\bR}\{c\1_{\{\alpha_{t}(V-c1_{\{t\}})\leq x\}}\},
\end{equation}
is a weakly acceptance (resp. weakly rejection) time consistent dynamic LM-measure.
\end{proposition}
In the proofs given in the Appendix, we will use the representation
\begin{equation}\label{eq:DAItoDCRM2}
\phi^{x}_{t}(V)(\omega)=\inf\{c\in\bR: \alpha_{t}(V-c1_{\{t\}})(\omega)\leq x\},
\end{equation}
rather than \eqref{eq:DAItoDCRM}, as it is more convenient.

This type of dual representations, i.e.  \eqref{eq:DCRMtoDAI} and \eqref{eq:DAItoDCRM}, first appeared in \cite{ChernyMadan2009}, where the authors studied static (one period of time) scale invariant measures. Subsequently, in \cite{BCZ2010}, the authors extended these results to the case of stochastic processes with special emphasis on time consistency property.
In contrast to \cite{BCZ2010}, we consider an arbitrary probability space, not just a finite one.

{We conclude this section by presenting two examples that illustrate the concept of semi-weak time consistency and show the connection between maps introduced in Propositions~\ref{prop:DCRMtoDAI}~and~\ref{prop:DAItoDCRM}. For more examples see \cite{BCP2014}.}

\begin{example}[Dynamic Gain Loss Ratio] \label{ex:4}
Dynamic Gain Loss Ratio (dGLR) is a popular measure of performance, which essentially improves on some drawbacks of Sharpe Ratio (such as penalizing for positive returns), and it is equal to the ratio of expected return over expected losses. Formally, for $\cX=\bV^{1}$, dGLR is defined as
\begin{equation}\label{e:dGLT}
\varphi_{t}(V):=\left\{
\begin{array}{ll}
\frac{E[\sum_{i=t}^{T}V_{i}|\cF_{t}]}{E[(\sum_{i=t}^{T}V_{i})^{-}|\cF_{t}]}, &\quad \textrm{if }  E[\sum_{i=t}^{T}V_{i}|\cF_{t}]>0,\\
0, &\quad \textrm{otherwise}.
\end{array}\right.
\end{equation}
For various properties and dual representations of dGLR see for instance \cite{BCZ2010,BCDK2013}.
In \cite{BCZ2010}, the authors showed that dGLR is both semi-weakly acceptance and semi-weakly rejection time consistent, although assuming that $\Omega$ is finite.
For sake of completeness we will show here that dGLR is semi-weakly acceptance time consistency; semi-weakly rejection time consistency is left to an interested reader as an exercise.

Assume that $t\in\bT\setminus\{T\}$,  and $V\in \bV^p$. In view of Proposition~\ref{prop:mart.ref}, it is enough to show that
\begin{equation}\label{eq:ex4.1}
\varphi_{t}(V)\geq \1_{\{V_{t}\geq 0\}}\Essinf_{t}(\varphi_{t+1}(V))+\1_{\{V_{t}< 0\}}(-\infty).
\end{equation}
On the set ${\{V_{t}< 0\}}$ the inequality \eqref{eq:ex4.1} is trivial.
Since $\varphi_{t}$ is non-negative and local, without loss of generality, we may assume that $\Essinf_{t}(\varphi_{t+1}(V))>0$.
Moreover, $\varphi_{t+1}(V)\geq\Essinf_{t}(\varphi_{t+1}(V))$, which implies
\begin{equation}\label{eq:ex4.2}
E[\sum_{i=t+1}^{T}V_{i}|\cF_{t+1}]\geq \Essinf_{t}(\varphi_{t+1}(V))\cdot E[(\sum_{i=t+1}^{T}V_{i})^{-}|\cF_{t+1}].
\end{equation}
Using \eqref{eq:ex4.2} we obtain
\begin{align}
\1_{\{V_{t}\geq 0\}}E[\sum_{i=t}^{T}V_{i}|\cF_{t}] & \geq \1_{\{V_{t}\geq 0\}}E[E[\sum_{i=t+1}^{T}V_{i}|\cF_{t+1}]|\cF_{t}]\nonumber\\
& \geq \1_{\{V_{t}\geq 0\}}\Essinf_{t}(\varphi_{t+1}(V))\cdot E[\1_{\{V_{t}\geq 0\}}E[(\sum_{i=t+1}^{T}V_{i})^{-}|\cF_{t+1}]|\cF_{t}]\nonumber\\
& \geq \1_{\{V_{t}\geq 0\}}\Essinf_{t}(\varphi_{t+1}(V))\cdot E[(\sum_{i=t}^{T}V_{i})^{-}|\cF_{t}].\label{eq:ex4.3}
\end{align}
{Note that $\Essinf_{t}(\varphi_{t+1}(V))>0$ implies that $\varphi_{t+1}(V)>0$, and thus $E[\sum_{i=t+1}^{T}V_{i}|\cF_{t+1}]>0$.
Hence, on the set $\set{V_t\geq0}$, we have
$$
\1_{\set{V_t\geq 1}}E[ \sum_{i=t}^T V_i |\cF_t ] \geq  \1_{\set{V_t\geq 1}} E[E[\sum_{i=t+1}^{T}V_{i}|\cF_{t+1}]|\cF_t]>0.
$$
Combining this and \eqref{eq:ex4.3}, we conclude the proof. }
\end{example}

\begin{example}[Dynamic RAROC for processes] \label{ex:3}
Risk Adjusted Return On Capital (RAROC) is a popular scale invariant measure of performance; see \cite{ChernyMadan2009} for a study of static RAROC, and  \cite{BCZ2010} for its extension to the dynamic setup. We consider the space $\cX=\bV^{1}$, and we fix  $\varepsilon\in (0,1)$.
Dynamic RAROC, at level $\varepsilon$, is the family $\{\alpha_{t}\}_{t\in\bT}$, with $\alpha_{t}$ given by
\begin{equation}\label{eq:dRAROC}
\alpha_{t}(V):=\left\{
\begin{array}{ll}
\frac{E[\sum_{i=t}^{T}V_{i}|\cF_{t}]}{-\rho^{\varepsilon}_{t}(V)}&\quad \textrm{if }  E[\sum_{i=t}^{T}V_{i}|\cF_{t}]>0,\\
0 &\quad \textrm{otherwise},
\end{array}\right.
\end{equation}
where $\rho^{\varepsilon}_{t}(V)=\essinf\limits_{Z\in \cD^{\varepsilon}_{t}}E[Z\sum_{i=t}^{T}V_{i}|\cF_{t}]$, and where the family of sets $\{D^{\varepsilon}_{t}\}_{t\in\bT}$ is defined by\footnote{The family $\{D^{\varepsilon}_{t}\}_{t\in\bT}$ represents risk scenarios, which define the dynamic version of the conditional value at risk at level $\varepsilon$ (cf. \cite{Cherny2010}).}
\begin{equation}\label{eq:det.sets.coh}
D^{\varepsilon}_{t}:=\{Z\in L^{1}: 0\leq Z\leq\varepsilon^{-1},\ E[Z|\cF_{t}]=1\}.
\end{equation}
We use the convention $\alpha_{t}(V)=+\infty$, if $\rho_{t}(V)\geq 0$.
In \cite{BCZ2010} it was shown that dynamic RAROC is a dynamic acceptability index for processes.
Moreover, it admits the following dual representation (cf. \eqref{eq:DCRMtoDAI2}): for any fixed $t\in\bT$,
$$
\alpha_{t}(V)=\sup\{x\in\bR_{+}: \phi^{x}_{t}(V)\geq 0 \},
$$
where $\phi^{x}_{t}(V)=\essinf\limits_{Z\in \cB^{x}_{t}}E[Z(\sum_{i=t}^{T}V_{i})|\cF_{t}]$, with
$$
\cB_{t}^{x}=\{Z\in L^{1}: Z=\frac{1}{1+x}+\frac{x}{1+x}Z_{1},\ \textrm{for some } Z_{1}\in \cD_{t}^{\varepsilon}\}.
$$
It is easy to check, that the family $\{\phi^{x}_{t}\}_{t\in\bT}$ is a dynamic coherent risk measure for processes, see \cite{BCZ2010} for details. Since $1\in \cD_{t}^{\varepsilon}$, we also get that $\{\phi_{t}^{x}\}_{t\in\bT}$ is increasing in $x\in\bR_{+}$.

Moreover, it is known that $\{\phi_{t}^{x}\}_{t\in\bT}$ is weakly acceptance time consistent but not weakly rejection time consistent, for any fixed $x\in\bR_{+}$ (see \cite[Example 1]{BCP2014}). Thus, using Propositions~\ref{prop:DCRMtoDAI}~and~\ref{prop:DAItoDCRM} we immediately conclude that $\{\phi^{x}_{t}\}_{t\in\bT}$ is semi-weakly acceptance time consistent and not semi-weakly rejection time consistent.

\end{example}

\begin{appendix}

\section{Appendix}

\subsection[Conditional expectation and essential supremum/infimum]{Conditional expectation and essential supremum/infimum on $\bar{L}^{0}$}\label{A:cond}
First, we will present some elementary properties of the generalized conditional expectation.
\begin{proposition}\label{pr:condexp}
For any $X,Y\in\bar{L}^{0}$ and $s,t\in\bT$, $s>t$ we get
\begin{enumerate}[1)]
\item $E[\lambda X|\cF_{t}]\leq \lambda E[X|\cF_{t}]$ for $\lambda\in L^{0}_{t}$, and $E[\lambda X|\cF_{t}]=\lambda E[X|\cF_{t}]$ for $\lambda\in L^{0}_{t}$, $\lambda\geq 0$;
\item $E[X|\cF_{t}]\leq E[E[X|\cF_{s}]|\cF_{t}]$, and $E[X|\cF_{t}]=E[E[X|\cF_{s}]|\cF_{t}]$ for $X\geq 0$;
\item $E[X|\cF_{t}]+E[Y|\cF_{t}]\leq E[X+Y|\cF_{t}]$, and $E[X|\cF_{t}]+E[Y|\cF_{t}]=E[X+Y|\cF_{t}]$ if $X,Y\geq 0$;
\end{enumerate}
\end{proposition}

\begin{remark}
All inequalities in Proposition~\ref{pr:condexp} can be strict.
Assume that $t=0$ and $k,s\in\bT$, $k>s > 0$, and let $\xi \in L^{0}_{k}$ be such that $\xi=\pm1$, $\xi$ is independent of $\cF_s$, and $P(\xi=1)=P(\xi=-1)=1/2$.
We consider $Z\in L_{s}^{0}$ such that $Z\geq 0$, and $E[Z]=\infty$.
By taking $\lambda=-1$, $X=\xi Z$ and $Y=-X$, we get strict inequalities in 1), 2) and 3).
\end{remark}

Next, we will discuss some important features of conditional essential infimum and conditional essential supremum, in the context of $\bar{L}^{0}$.

Before that, we will recall the definition of conditional essential infimum for bounded random variables. For $X\in L^{\infty}$ and $t\in\bT$, we will denote by $\Essinf_{t}X$ the unique (up to a set of probability zero), $\cF_{t}$-measurable random variable, such that for any $A\in\cF_{t}$, the following equality holds true
\begin{equation}\label{eq:essinf.b}
\essinf_{\omega\in A}X=\essinf_{\omega\in A}(\Essinf_{t}X).
\end{equation}
We will call this random variable the \textit{$\cF_{t}$-conditional essential infimum of $X$}. We refer the reader to \cite{BarronCardaliaguetJensen2003} for a detailed proof of the existence and uniqueness of the conditional essential infimum.
We will call $\Esssup_{t}(X):=-\Essinf_{t}(-X)$ the  \textit{$\cF_{t}$-conditional essential supremum of $X\in L^\infty$}.

As stated in the preliminaries, we extend these two notions to the space $\bar{L}^0$.
 For any $t\in\bT$ and $X\in\bar{L}^{0}$, we define the $\cF_{t}$-conditional essential infimum by
\begin{equation}\label{eq:A22}
\Essinf_{t}X:=\lim_{n\to\infty}\Big[\Essinf_{t}(X^{+}\wedge n)\Big]-\lim_{n\to\infty}\Big[\Esssup_{t}(X^{-}\wedge n)\Big],
\end{equation}
and respectively we put $\Esssup_{t}(X):=-\Essinf_{t}(-X)$.

\begin{remark}
Extending the function $\arctan$ to $[-\infty,\infty]$ by continuity, and observing that $\arctan X\in L^\infty$ for any $X\in \bar{L}^{0}$, one can naturally extend conditional essential infimum to $ \bar{L}^{0}$ by setting
\[
\Essinf_{t}X= \arctan^{-1}[\Essinf_{t}(\arctan X)].
\]
\end{remark}
We proceed with the following result:
\begin{proposition}\label{pr:essinf}
For any $X,Y\in \bar{L}^{0}$, $s,t \in \bT, s\geq t$, and $A\in\cF_{t}$ we have
\begin{enumerate}[1)]
\item $\essinf_{\omega\in A}X=\essinf_{\omega\in A}(\Essinf_{t}X)$;
\item If $\essinf_{\omega\in A}X=\essinf_{\omega\in A}U$ for some $U\in \bar{L}^{0}_t$, then $U=\Essinf_{t}X$;
\item $X\geq \Essinf_{t}X$;
\item  If $Z\in\bar{L}^{0}_{t}$, is such that $X\geq Z$, then $\Essinf_{t}X\geq Z$;
\item If $X\geq Y$, then $\Essinf_{t}X\geq \Essinf_{t}Y$;
\item $\1_{A}\Essinf_{t}X=\1_{A}\Essinf_{t}(\1_{A}X)$;
\item $\essinf_{s}X\geq \Essinf_{t}X$;
\end{enumerate}
The analogous results are true for $\{\Esssup_{t}\}_{t\in\bT}$.
\end{proposition}
The proof for the case $X,Y\in L^{\infty}$ can be found in~\cite{BarronCardaliaguetJensen2003}. Since for any $n\in\bN$ and $X,Y\in \bar{L}^{0}$ we get $X^{+}\wedge n\in L^{\infty}$, $X^{-}\wedge n\in L^{\infty}$ and $X^{+}\wedge X^{-}=0$, the extension of the proof to the case $X,Y\in \bar{L}^{0}$ is straightforward, and we omit it here.

\begin{remark}\label{rem:essinf.def2}  Similarly to~\cite{BarronCardaliaguetJensen2003},  the conditional essential infimum $\essinf_{t}(X)$  can be alternatively defined as the largest $\cF_{t}$-measurable random variable, which is smaller than $X$, i.e. properties 3) and 4) from Proposition~\ref{pr:essinf} are characteristic properties for conditional essential infimum.
\end{remark}

Next, we define the generalized versions of $\essinf$ and $\esssup$ of a (possibly uncountable) family of random variables:
For $\{X_{i}\}_{i\in I}$, where $X_{i}\in\bar{L}^{0}$, we let
\begin{equation}\label{eq:cond222}
\essinf_{i\in I}X_{i}:=\lim_{n\to\infty}\Big[\Essinf_{i\in I}(X_{i}^{+}\wedge n)\Big]-\lim_{n\to\infty}\Big[\Esssup_{i\in I}(X_{i}^{-}\wedge n)\Big].
\end{equation}
Note that, in view of \cite[Appendix~A]{KaratzasShreve1998},  $\Essinf_{i\in I}X_{i}\wedge n $ and $\Esssup_{i\in I}X_{i}\wedge n $ are well defined, so that  $\essinf_{i\in I}X_{i}$ is well defined. It needs to be observed that the operations of the right hand side of \eqref{eq:cond222} preserve measurability. In particular, if $X_{i}\in\cF_{t}$ for all $i\in I$, then $\essinf_{i\in I}X_{i}\in\cF_{t}$.

Furthermore, if for any $i,j\in I$, there exists $k\in I$, such that $X_{k}\leq X_{i}\wedge X_{j}$, then there exists a sequence $i_n\in I, n\in\bN$, such that $\{X_{i_n}\}_{n\in\bN}$ is nonincreasing and $\essinf_{i\in I}X_{i}=\inf_{n\in\bN}X_{i_n}=\lim_{n\to\infty}X_{i_n}$.
Analogous results hold true for $\esssup_{i\in I}X_{i}$.

\subsection{Proofs}\label{A:proofs}

\subsubsection*{Proof of Proposition~\ref{prop:mart.ref}.}\label{prop:mart.ref.a}
\begin{proof}
Let $\mu$ be an update rule.

\smallskip
\noindent 1) The implication  ($\Rightarrow$) follows immediately, by taking in the definition of acceptance time consistency $m_s=\varphi_s(X)$.

\smallskip
\noindent ($\Leftarrow$) Assume that $\varphi_{t}(X)\geq \mu_{t,s}(\varphi_{s}(X),X)$,  for any $s,t\in\bT, s>t$, and $X\in\cX$.
Let $m_{s}\in\bar{L}^{0}_{s}$ be such that $\varphi_{s}(X)\geq m_{s}$. Using monotonicity of $\mu$, we get $\varphi_{t}(X)\geq\mu_{t,s}(\varphi_{s}(X),X)\geq \mu_{t,s}(m_{s},X)$.

\smallskip
\noindent 2) The proof is similar to 1).
\end{proof}

\subsubsection*{Proof of Proposition~\ref{prop:benchmark}.}\label{prop:benchmark.a}

\begin{proof}
We do the proof only for acceptance time consistency. The proof for rejection time consistency is analogous.

\smallskip\noindent
\textit{Step 1.} We will show that $\varphi$ is acceptance time consistent with respect to $\cY$, if and only if
\begin{equation}\label{eq:bench3}
\1_A\varphi_{s}(X)\geq \1_A\varphi_{s}(Y)\quad \Longrightarrow\quad \1_A\varphi_{t}(X)\geq \1_A\varphi_{t}(Y),
\end{equation}
for all $s\geq t$, $X\in L^p$, $Y\in \cY_{s}$ and $A\in\cF_{t}$. For sufficiency it is enough to take $A=\Omega$. For necessity let us assume that
\begin{equation}\label{eq:bench.pr1}
\1_A\varphi_s(X)\geq \1_A\varphi_s(Y).
\end{equation}
Using locality of $\varphi$, we get that \eqref{eq:bench.pr1} is equivalent to
\[
\1_A\varphi_s(\1_A X+\1_{A^c}Y)+\1_{A^c}\varphi_s(\1_A X+\1_{A^c}Y)\geq \1_A\varphi_s(Y)+\1_{A^c}\varphi_s(Y),
\]
and consequently to $\varphi_s(\1_A X+\1_{A^c}Y)\geq \varphi_s(Y)$.
Thus, using \eqref{eq:benchmark2}, we get
\[
\varphi_s(\1_A X+\1_{A^c}Y)\geq \varphi_s(Y) \quad \Longrightarrow\quad \varphi_t(\1_A X+\1_{A^c}Y)\geq \varphi_t(Y).
\]
By the same arguments we get that $\varphi_t(\1_A X+\1_{A^c}Y)\geq \varphi_t(Y)$ is equivalent to $\1_A\varphi_t(X)\geq \1_A\varphi_t(Y)$, which concludes this part of the proof.

\smallskip\noindent
{\it Step 2.}  Now we demonstrate that $\varphi$ is acceptance time consistent with respect to $\cY$ if and only if $\varphi$ is acceptance time consistent with respect to the family $\widehat \cY=\{\widehat \cY_t\}_{t\in\bT}$ of benchmark sets given by
\begin{equation}\label{eq:bench.local}
\widehat{\cY}_{t}:=\{\1_{A}Y_1+\1_{A^c}Y_2\,:\, Y_1,Y_2\in\cY_{t},\, A\in\cF_t\}.
\end{equation}
Noting that for any $t\in\bT$ we have $\cY_t\subseteq \widehat \cY_t$, we get the sufficiency part. For necessity let us assume that
\begin{equation}\label{eq:bench.pr2}
\varphi_s(X)\geq \varphi_s(Y),
\end{equation}
for some $Y\in\widehat \cY_t$. In view of \eqref{eq:bench.local} we conclude that there exists $A\in\cF_t$ and $Y_1,Y_2\in\cY_s$, such that $Y=\1_{A}Y_1+\1_{A^c}Y_2$. Consequently, using locality of $\varphi$, and the fact that \eqref{eq:bench.pr2} is equivalent to
\[
\1_A\varphi_s(X)+\1_{A^c}\varphi_s(X)\geq \1_{A}\varphi_s(\1_{A}Y_1+\1_{A^c}Y_2)+\1_{A^c}\varphi_s(\1_{A}Y_1+\1_{A^c}Y_2),
\]
we deduce that \eqref{eq:bench.pr2} is equivalent to
\[
\1_A\varphi_s(X)+\1_{A^c}\varphi_s(X)\geq \1_{A}\varphi_s(Y_1)+\1_{A^c}\varphi_s(Y_2).
\]
As the sets $A$ and $A^c$ are disjoint, using \eqref{eq:bench3} twice, we get
\[
\1_A\varphi_t(X)+\1_{A^c}\varphi_t(X)\geq \1_{A}\varphi_t(Y_1)+\1_{A^c}\varphi_t(Y_2).
\]
By similar arguments as before, we get that the above inequality is equivalent to $\varphi_t(X)\geq \varphi_t(Y)$, that  concludes this part of the proof.

\smallskip\noindent
\textit{Step 3.} For any $m_s\in \bar{L}^0_s$, we set
\[
\mu_{t,s}(m_s):=\esssup_{A\in\cF_{t}}\Big[\1_{A}\esssup_{Y\in \cY^{-}_{A,s}(m_s)}\varphi_{t}(Y)+\1_{A^{c}}(-\infty)\Big],
\]
where $\cY^{-}_{A,s}(m_s):=\{Y\in\widehat{\cY}_s: \1_A m_s\geq \1_A\varphi_s(Y)\}$, and show that the corresponding family of maps $\mu$ is  a projective update rule.

\smallskip
\noindent \textit{Adaptiveness.} For any $m_s\in\bar{L}^0_s$, $\esssup$ of the set of $\cF_{t}$-measurable random variables $\{\varphi_{t}(Y)\}_{Y\in \cY^{-}_{A,s}(m_s)}$ is $\cF_{t}$-measurable (see~\cite{KaratzasShreve1998}, Appendix A), which implies that $\mu_{t,s}(m_s)\in\bar{L}^{0}_{t}$.

\smallskip
\noindent \textit{Monotonicity.} If $m_s\geq m_s'$, then for any $A\in\cF_{t}$ we get $\cY^{-}_{A,s}(m_s)\supseteq \cY^{-}_{A,s}(m_s')$,
which implies $\mu_{t,s}(m_s)\geq\mu_{t,s}(m_s')$.

\smallskip
\noindent \textit{Locality.} Let $B\in\cF_{t}$, and $m_s\in\bar{L}^0_s$.
It is enough to consider $A\in\cF_{t}$, such that $\cY^{-}_{A,s}(m_s)\neq\emptyset$, as otherwise we get
\[
\Big[\1_{A}\esssup_{Y\in \cY^{-}_{A,s}(m_s)}\varphi_{t}(Y)+\1_{A^{c}}(-\infty)\Big]\equiv -\infty.
\]
For any such $A\in\cF_{t}$, we get
\begin{equation}\label{eq:ppp1pp1}
\1_{A\cap B}\esssup_{Y\in  \cY^{-}_{A,s}(m_s)}\varphi_{t}(Y)=\1_{A\cap B}\esssup_{Y\in  \cY^{-}_{A\cap B,s}(m_s)}\varphi_{t}(Y).
\end{equation}
Indeed, since $ \cY^{-}_{A,s}(m_s)\subseteq  \cY^{-}_{A\cap B,s}(m_s)$, we have
$$
\1_{A\cap B}\esssup_{Y\in  \cY^{-}_{A,s}(m_s)}\varphi_{t}(Y)\leq \1_{A\cap B}\esssup_{Y\in  \cY^{-}_{A\cap B,s}(m_s)}\varphi_{t}(Y).
$$
On the other hand, for any $Y\in \cY^{-}_{A\cap B,s}(m_s)$, and for a fixed $Z\in \cY^{-}_{A,s}(m_s)$,   in view of \eqref{eq:bench.local}, we obtain
$$
\1_{B}Y+\1_{B^{c}}Z\in  \cY^{-}_{A,s}(m_s).
$$
Thus, using locality of $\varphi_{t}$, we deduce
$$
\1_{A\cap B}\esssup_{Y\in \cY^{-}_{A\cap B,s}(m_s)}\varphi_{t}(Y)= \1_{A\cap B}\esssup_{Y\in  \cY^{-}_{A\cap B,s}(m_s)}\1_{B}\varphi_{t}(\1_{B}Y+\1_{B^{c}}Z)\leq\1_{A\cap B}\esssup_{Y\in  \cY^{-}_{A,s}(m_s)}\varphi_{t}(Y),
$$
which proves \eqref{eq:ppp1pp1}.
Now, note that $ \cY^{-}_{A\cap B,s}(m_s)= \cY^{-}_{A\cap B,s}(\1_B m_s)$, and thus
\begin{equation}\label{eq:ppp1pp2}
\1_{A}\esssup_{Y\in \cY^{-}_{A\cap B,s}(m_s)}\varphi_{t}(Y)=\1_{A}\esssup_{Y\in  \cY^{-}_{A\cap B,s}(\1_B m_s)}\varphi_{t}(Y).
\end{equation}
Combining \eqref{eq:ppp1pp1}, \eqref{eq:ppp1pp2}, and the fact that $ \cY^{-}_{A,s}(m_s)\neq\emptyset$ implies $ \cY^{-}_{A,s}(\1_B m_s)\neq\emptyset$, we obtain the following chain of equalities
\begin{align*}
\1_{B}\mu_{t,s}(m_s) &= \1_{B}\esssup_{A\in\cF_{t}}\Big[\1_{A}\esssup_{Y\in  \cY^{-}_{A,s}(m_s)}\varphi_{t}(Y)+\1_{A^{c}}(-\infty)\Big]\\
& = \1_{B}\esssup_{A\in\cF_{t}}\Big[\1_{A\cap B}\esssup_{Y\in  \cY^{-}_{A,s}(m_s)}\varphi_{t}(Y)+\1_{A^{c}\cap B}(-\infty)\Big]\\
& = \1_{B}\esssup_{A\in\cF_{t}}\Big[\1_{A\cap B}\esssup_{Y\in  \cY^{-}_{A\cap B,s}(m_s)}\varphi_{t}(Y)+\1_{A^{c}\cap B}(-\infty)\Big]\\
& = \1_{B}\esssup_{A\in\cF_{t}}\Big[\1_{A\cap B}\esssup_{Y\in  \cY^{-}_{A\cap B,s}(\1_B m_s)}\varphi_{t}(Y)+\1_{A^{c}\cap B}(-\infty)\Big]\\
& =\1_{B}\esssup_{A\in\cF_{t}}\Big[\1_{A}\esssup_{Y\in \cY^{-}_{A,s}(\1_B m_s)}\varphi_{t}(Y)+ \1_{A^{c}}(-\infty)\Big]\\
& = \1_{B}\mu_{t,s}(\1_B m_s).
\end{align*}
Thus, $\mu$ is an $X$-invariant update rule.

\smallskip\noindent
\textit{Step 4.} By locality of $\varphi$ and \eqref{eq:bench3}, we note that acceptance time consistency with respect to $\cY$ is equivalent to
\begin{equation}\label{eq:bench.pr5}
\varphi_{t}(X)\geq \esssup_{A\in\cF_{t}}\Big[\1_{A}\esssup_{Y\in \cY^{-}_{A,s}(\varphi_s(X))}\varphi_{t}(Y)+\1_{A^{c}}(-\infty)\Big].
\end{equation}
Thus, using \eqref{eq:accepTimeConsAlt}, we deduce that $\varphi$ satisfies \eqref{eq:benchmark2} if and only if $\varphi$ is time consistent with respect to the update rule $\mu$. Since \eqref{eq:accepTimeConsAlt} is  equivalent to \eqref{eq:bench.pr5}, we conclude the proof.
\end{proof}

\subsubsection*{Proof of Proposition~\ref{th:essinf}.}\label{th:essinf.a}
\begin{proof}
Monotonicity and locality of $\mu^{\inf}$ is a straightforward implication of Proposition~\ref{pr:essinf}. Thus, $\mu^{\inf}$ is an $sX$-invariant update rule. The projectivity comes straight from the definition (see Remark~\ref{rem:essinf.def2}).
\end{proof}

\subsubsection*{Proof of Proposition~\ref{pr:weak}.}\label{pr:weak.a}
 \begin{proof}
We will only show the proof for acceptance consistency.
The proof for rejection consistency is similar.
Let $\varphi$ be a dynamic LM-measure.


 \smallskip
 \noindent $1)\Rightarrow 2)$. Assume that $\varphi$ is weakly acceptance consistent, and  let $m_{t}\in \bar{L}^{0}_{t}$ be such that $\varphi_{s}(X)\geq m_{t}$.
 Then, using Proposition~\ref{prop:mart.ref},  we get $\varphi_{t}(X)\geq \Essinf_{t}(\varphi_{s}(X))\geq \Essinf_{t}(m_{t})=m_{t}$, and hence 2) is proved.

 \smallskip
 \noindent $2)\Rightarrow 1)$. By the definition of conditional essential infimum, $\Essinf_{t}(\varphi_{s}(X))\in\bar{L}^{0}_{t}$, for any $X\in L^p$, and $t,s\in\cT$.
  Moreover, by Proposition~\ref{pr:essinf}.(3), we have that  $\varphi_{s}(X)\geq\Essinf_{t}(\varphi_{s}(X))$.
  Using  2) with  $m_t=\Essinf_{t}(\varphi_{s}(X))$,  we immediately obtain $\varphi_{t}(X)\geq\Essinf_{t}(\varphi_{s}(X))$. Due to Proposition~\ref{prop:mart.ref}, this concludes the proof.

 \smallskip
 \noindent $2)\Leftrightarrow 3)$. Clearly 2) $\Rightarrow $ 3). Thus, it remains to show the converse implication. Since $\varphi$ is a monetary utility measure, then invoking locality of $\varphi$, we conclude that for any $m_{t}\in \bar{L}^{0}_{t}$, such that $\varphi_s(X)\geq m_t$, and for any $n\in\bN$, we have
\[
\varphi_s(\1_{\{m_{t}\in (-n,n)\}}(X-m_t))\geq 0.
\]
Now, in view of 3), we get that $\varphi_t(\1_{\{m_{t}\in (-n,n)\}}(X-m_t))\geq0$, and consequently
\[
\1_{\{m_{t}\in (-n,n)\}}\varphi_t(X)\geq\1_{\{m_{t}\in (-n,n)\}}m_{t}.
\]
Thus, 2) is proved on the $\cF_{t}$-measurable set $\{m_{t}\in(-\infty,\infty)\}=\bigcup_{n\in\bN}\{m_t\in (-n,n)\}$. On the set $\{m_t=-\infty\}$ the inequality $\varphi_{t}(X)\geq m_{t}$ is trivial. Finally, on the set $\{m_t=\infty\}$, in view of  the monotonicity of $\varphi$, we have that $\varphi_{s}(X)=\varphi_{t}(X)=\infty$, which implies 2). This concludes the proof.

 \end{proof}

\subsubsection*{Proof of Proposition~\ref{pr:essinf.cond.aa}.}\label{pr:essinf.cond.aa.a}

\begin{proof}
Let a family $\mu=\{\mu_{t}\}_{t\in\bT}$ of maps $\mu_{t}:\bar{L}^{0}\to \bar{L}^{0}_{t}$ be given by
\begin{equation}\label{eq:essinff}
\mu_{t}(m)=\essinf_{Z\in P_{t}}E[Zm|\cF_{t}]
\end{equation}
Before proving \eqref{eq:essinf.eq.pt}, we will need to prove some auxiliary facts about $\mu$.

First, let us show that $\mu$ is local and monotone. Let $t\in\bT$. Monotonicity is straightforward. Indeed,  let $m,m'\in\bar{L}^{0}$ be such that $m\geq m'$. For any $Z\in P_{t}$, using the fact that $Z\geq 0$, we get $Zm\geq Zm'$. Thus, $E[Zm|\cF_{t}]\geq E[Zm'|\cF_{t}]$ and consequently $\essinf_{Z\in P_{t}}E[Zm|\cF_{t}]\geq \essinf_{Z\in P_{t}}E[Zm'|\cF_{t}]$.
Next, for any $A\in\cF_{t}$ and $m\in\bar{L}^{0}$, by invoking  Proposition~\ref{pr:condexp}, convention $0\cdot\pm\infty=0$, and the fact that for any $Z_{1},Z_{2}\in P_{t}$ we have $\1_{A}Z_{1}+\1_{A^{c}}Z_{2}\in P_{t}$, we get
\begin{align*}
\1_{A}\mu_{t}(m) & = \1_{A}\essinf_{Z\in P_{t}}E[Zm|\cF_{t}]\\
& =\1_{A}\essinf_{Z\in P_{t}}(E[(\1_{A}Z)m|\cF_{t}]+E[(\1_{A^{c}}Z)m|\cF_{t}])\\
&=\1_{A}\essinf_{Z\in P_{t}}E[(\1_{A}Z)m|\cF_{t}]+\1_{A}\essinf_{Z\in P_{t}}E[(\1_{A^{c}}Z)m|\cF_{t}]\\
&=\1_{A}\essinf_{Z\in P_{t}}E[Z(\1_{A}m)|\cF_{t}]+\1_{A}\essinf_{Z\in P_{t}}\1_{A^{c}}E[Zm|\cF_{t}]\\
&=\1_{A}\mu_{t}(\1_{A}m),
\end{align*}
which proves locality.

Secondly, let us prove that
\begin{equation}\label{eq:A0}
m \geq \mu_{t}(m),
\end{equation}for any $m\in \bar{L}^{0}$. Let $m\in L^{0}$. For $\varepsilon\in (0,1)$ let\footnote{In the risk measure framework, it might be seen as the risk minimazing scenario for conditional $CV@R_{\varepsilon}$.}
\begin{equation}\label{eq:Zalpha}
Z_{\varepsilon}:=\1_{\set{m\leq q^{+}_{t}(\varepsilon)}}E[\1_{\set{m\leq q^{+}_{t}(\varepsilon)}}|\cF_{t}]^{-1}.
\end{equation}
where $q^{+}_{t}(\varepsilon)$ is  $\cF_{t}$-conditional (upper) $\varepsilon$ quantile of $m$, defined as
$$
q^+_{t}(\varepsilon):= \esssup\set{ Y\in L_{t}^0 \mid E[\1_{\set{m\leq Y}}|\cF_{t}]\leq\varepsilon}.
$$
For $\varepsilon\in (0,1)$, noticing that $Z_{\varepsilon}<\infty$, due to convention $0\cdot\infty=0$ and the fact that
$$
\{E[\1_{\set{m\leq q^{+}_{t}(\varepsilon)}}|\cF_{t}]=0\}\subseteq \{\1_{\set{m\leq q^{+}_{t}(\varepsilon)}}= 0\}\cup B,$$
for some $B$, such that $P[B]=0$, we conclude that $Z_{\varepsilon}\in P_{t}$.

Moreover, by the definition of $q^+_t(\varepsilon)$, there exists a sequence $Y_n\in L_{t}^0$, such that $Y_n\nearrow q^+_t(\varepsilon)$, and
$$
E[\1_{\set{m<Y_n}} \ |\ \cF_t]\leq \varepsilon.
$$
Consequently, by monotone convergence theorem, we have
$$
E[\1_{\set{m<q^+_t(\varepsilon)}} \ | \ \cF_t] \leq \varepsilon.
$$
Hence, we deduce
$$P[m<q^{+}_{t}(\varepsilon)]=E[\1_{\set{m<q^{+}_{t}(\varepsilon)}}]\leq E[E[\1_{\set{m<q^{+}_{t}(\varepsilon)}}|\cF_{t}]]\leq E[\varepsilon]=\varepsilon,$$
which implies that
\begin{equation}\label{eq:A6}
P[m\geq q^{+}_{t}(\varepsilon)]\geq (1-\varepsilon).
\end{equation}
On the other hand
\begin{align*}
\1_{\set{m\geq q^{+}_{t}(\varepsilon)}}m & \geq \1_{\set{m\geq q^{+}_{t}(\varepsilon)}} q_{t}^{+}(\varepsilon)= \1_{\set{m\geq q^{+}_{t}(\varepsilon)}}  q_{t}^{+}(\alpha)E[Z_{\varepsilon}|\cF_{t}]\\
& \geq \1_{\set{m\geq q^{+}_{t}(\varepsilon)}} E[Z_{\varepsilon} q_{t}^{+}(\varepsilon)|\cF_{t}]\geq \1_{\set{m\geq q^{+}_{t}(\varepsilon)}} E[Z_{\varepsilon}m|\cF_{t}],
\end{align*}
which combined with \eqref{eq:A6}, implies that
\begin{equation}\label{eq:1alpha}
P\Big[m\geq E[Z_{\varepsilon}m|\cF_{t}]\Big]\geq 1-\varepsilon.
\end{equation}
Hence, using \eqref{eq:1alpha}, and the fact that
$$
E[Z_{\varepsilon}m|\cF_{t}]\geq\mu_{t}(m), \quad \varepsilon\in(0,1),
$$
we get that
$$
P[m\geq \mu_{t}(m)]\geq 1-\varepsilon.
$$
Letting $\varepsilon\to 0$, we conclude that \eqref{eq:A0} holds true for $m\in L^{0}$.

Now, assume that $m\in\bar{L}^{0}$, and let $A:=\{E[\1_{\{m=-\infty\}}|\cF_{t}]=0\}$.
Similar to the arguments above, we get
$$
\1_{A}m\geq \mu_{t}(\1_{A}m).
$$
Since $\mu_{t}(0)=0$, and due to locality of $\mu_{t}$, we deduce
\begin{equation}\label{eq:A10}
\1_{A}m\geq \mu_{t}(\1_{A}m)=\1_{A}\mu_{t}(\1_{A}m)=\1_{A}\mu_{t}(m).
\end{equation}
Moreover, taking $Z=1$ in \eqref{eq:essinff}, we get
\begin{equation}\label{eq:A11}
\1_{A^{c}}m\geq \1_{A^{c}}(-\infty)= \1_{A^{c}}E[m|\cF_{t}]\geq \1_{A^{c}}\mu_{t}(m).
\end{equation}
Combining \eqref{eq:A10} and \eqref{eq:A11}, we conclude the proof of \eqref{eq:A0}  for all $m\in\bar{L}^{0}$.

Finally, we will show that $\mu$ defined as in \eqref{eq:essinff} satisfies property 1) from Proposition~\ref{pr:essinf}, which will consequently imply equality \eqref{eq:essinf.eq.pt}. Let $m\in\bar{L}^{0}$ and $A\in\cF_{t}$. From the fact that $m\geq \mu_{t}(m)$, we get
$$\essinf_{\omega\in A}m\geq\essinf_{\omega\in A}\mu_{t}(m).$$
On the other hand, we know that $\1_{A}\essinf_{\omega\in A}m\leq \1_{A}m$  and $\1_{A}\essinf_{\omega\in A}m\in\bar{L}^{0}_{t}$, so
\begin{align*}
\essinf_{\omega\in A}m & =\essinf_{\omega\in A}(\1_{A}\essinf_{\omega\in A}m)=\essinf_{\omega\in A}(\1_{A}\mu_{t}(\1_{A}\essinf_{\omega\in A}m))\leq\\
& \leq \essinf_{\omega\in A}(\1_{A}\mu_{t}(\1_{A}m))=\essinf_{\omega\in A}(\1_{A}\mu_{t}(m))=\essinf_{\omega\in A}\mu_{t}(m)
\end{align*}
which proves the equality. The proof for $\Esssup_{t}$ is similar and we omit it here.
This concludes the proof.

\end{proof}

\subsubsection*{Proof of Proposition~\ref{pr:UDMprop}.}\label{pr:UDMprop.a}
\begin{proof}
Using Proposition~\ref{pr:essinf}, for any $t,s\in\mathbb{T}$, $s>t$, and any $X\in  L^p$, we get
$$
\varphi_{t}(X)\geq \mu_{t}(\varphi_{s}(X))\geq \mu_{t}(\Essinf_{s}(\varphi_{s}(X)))\geq\mu_{t}(\Essinf_{t}(\varphi_{s}(X)))=\Essinf_{t}(\varphi_{s}(X)).
$$
The proof for rejection time consistency is similar.
\end{proof}

\subsubsection*{Proof of Proposition~\ref{pr:proc.weak.semi}.}\label{pr:proc.weak.semi.a}

 \begin{proof}
We will only show the proof for acceptance consistency.
The proof for rejection consistency is similar.
Let $\varphi$ be a dynamic LM-measure.

 \smallskip
 \noindent $1)\Leftrightarrow 2)$. This is a direct implication of Proposition~\ref{prop:mart.ref}.

 \smallskip
 \noindent $2)\Rightarrow 3)$. Assume that $\varphi$ is semi-weakly acceptance consistent. Let $V\in\bV^p$ and $m_{t}\in \bar{L}^{0}_{t}$ be such that $\varphi_{t+1}(V)\geq m_{t}$ and $V_{t}\geq 0$.  Then, by monotonicity of $\mu^{\inf}_{t}$, we have
 $$
 \varphi_{t}(V)\geq \1_{\set{V_{t}\geq 0}}\mu^{\inf}_{t}(\varphi_{t+1}(V))\geq \mu^{\inf}_{t}(m_{t})= \Essinf_{t}(m_{t})=m_{t},
 $$
 and hence 3) is proved.

 \smallskip
 \noindent $3)\Rightarrow 2)$. Let $V\in\bV^p$. We need to show that
 \begin{equation}\label{eq:semi.A}
 \varphi_{t}(V)\geq 1_{\{V_{t}\geq 0\}}\mu^{\inf}_{t}(\varphi_{t+1}(V))+1_{\{V_{t}< 0\}}(-\infty).
\end{equation}
 On the set $\{V_{t}<0\}$ inequality \eqref{eq:semi.A} is trivial. We know that
 $$
 (\1_{\{V_{t}\geq 0\}}\cdot_{t}V)_{t}\geq 0\quad \textrm{and}\quad \varphi_{t+1}(\1_{\{V_{t}\geq 0\}}\cdot_{t}V)\geq \Essinf_{t}\varphi_{t+1}(\1_{\{V_{t}\geq 0\}}\cdot_{t}V).
 $$
 Thus, for $m_{t}=\Essinf_{t}\varphi_{t+1}(\1_{\{V_{t}\geq 0\}}\cdot_{t}V)$, using locality of $\varphi$ and $\mu^{\inf}$ as well as 3), we get
$$
\1_{\{V_{t}\geq 0\}}\varphi_{t}(V)=\1_{\{V_{t}\geq 0\}}\varphi_{t}(\1_{\{V_{t}\geq 0\}}\cdot_{t}V)\geq \1_{\{V_{t}\geq 0\}} m_{t}=\1_{\{V_{t}\geq 0\}}\mu_{t}^{\inf}(\varphi_{t+1}(V)),
$$
and hence \eqref{eq:semi.A} is proved on the set $\{V_{t}\geq 0\}$. This conclude the proof of 2).
 \end{proof}

\subsubsection*{Proof of Proposition~\ref{prop:DCRMtoDAI}}
\begin{proof}
The proof of locality and monotonicity of \eqref{eq:DCRMtoDAI} is straightforward (see \cite{BCZ2010} for details). Let us assume that $\{\phi^{x}_{t}\}_{t\in\bT}$ is weakly acceptance time consistent. Using counterpart of Proposition~\ref{pr:weak} for stochastic processes (see~\cite{BCP2014}) we get
\begin{align*}
1_{\{V_{t}\geq 0\}}\alpha_{t}(V) & = 1_{\{V_{t}\geq 0\}}\Big(\sup\{x\in\bR_{+}: 1_{\{V_{t}\geq 0\}}\phi_{t}^{x}(V)\geq0\}\Big)\\
& \geq 1_{\{V_{t}\geq 0\}}\Big(\sup\{x\in\bR_{+}: 1_{\{V_{t}\geq 0\}}[\Essinf_{t}\phi_{t+1}^{x}(V)+V_{t}]\geq0\}\Big)\\
& \geq 1_{\{V_{t}\geq 0\}}\Big(\sup\{x\in\bR_{+}: 1_{\{V_{t}\geq 0\}}\Essinf_{t}\phi_{t+1}^{x}(V)\geq0\}\Big)\\
& = 1_{\{V_{t}\geq 0\}}\Essinf_{t}\Big(\sup\{x\in\bR_{+}: 1_{\{V_{t}\geq 0\}}\phi_{t+1}^{x}(V)\geq0\}\Big)\\
& = 1_{\{V_{t}\geq 0\}}\Essinf_{t}\alpha_{t+1}(V).
\end{align*}
This leads to
$$\alpha_{t}(V)\geq 1_{\{V_{t}\geq 0\}}\Essinf_{t}\alpha_{t+1}(V) +1_{\{V_{t}< 0\}}(-\infty),$$
which, by Proposition~\ref{pr:proc.weak.semi}, is equivalent to the semi-weak rejection time consistency. The proof of the weak acceptance time consistency is similar.
\end{proof}

\subsubsection*{Proof of Proposition~\ref{prop:DAItoDCRM}}
\begin{proof}[Proof]

The proof of locality and monotonicity of \eqref{eq:DAItoDCRM} is straightforward (see \cite{BCZ2010} for details). Let us prove the weak acceptance time consistency. Assume that $\{\alpha_{t}\}_{t\in\bT}$ is semi-weakly acceptance time consistent. Using Proposition~\ref{prop:mart.ref} we get
\begin{align*}
\phi^{x}_{t}(V) & = \inf\{c\in\bR: \alpha_{t}(V-c1_{\{t\}})\leq x\}\\
& = \inf\{c\in\bR: \alpha_{t}(V-c1_{\{t+1\}})\leq x\}\\
& = \inf\{c\in\bR: \alpha_{t}(V-c1_{\{t+1\}}-V_{t}1_{\{t\}})\leq x\}+V_{t}\\
& \geq \inf\{c\in\bR: 1_{\{0\geq 0\}}\Essinf_{t}\alpha_{t+1}(V-c1_{\{t+1\}}-V_{t}1_{\{t\}}) +1_{\{0< 0\}}(-\infty)\leq x\}+V_{t}\\
& = \inf\{c\in\bR: \Essinf_{t}\alpha_{t+1}(V-c1_{\{t+1\}})\leq x\}+V_{t}\\
& = \Essinf_{t}\big(\inf\{c\in\bR:\alpha_{t+1}(V-c1_{\{t+1\}})\leq x\}\big)+V_{t}\\
& = \Essinf_{t}\phi_{t+1}^{x}(V)+V_{t},
\end{align*}
which,  is equivalent to the weak acceptance  time consistency of $\phi$. The proof of the rejection time consistency is similar.
\end{proof}

\subsubsection*{Proof of Proposition~\ref{pr:condexp}.}\label{pr:condexp.a}
\begin{proof}
First note that for any $X,Y\in \bar{L}^{0}$, $\lambda\in L^{0}_{t}$, such that $X,Y, \lambda \geq 0$, and for any $s,t\in\bT$, $s>t$, by Monotone Convergence Theorem, and using the convention $0\cdot\pm\infty=0$,  we get
\begin{eqnarray}
E[\lambda X|\cF_{t}]& = &\lambda E[X|\cF_{t}];\label{eq:1.1}\\
E[X|\cF_{t}]& = &E[E[X|\cF_{s}]|\cF_{t}];\label{eq:1.2}\\
E[X|\cF_{t}]+E[Y|\cF_{t}]&=&E[X+Y|\cF_{t}].\label{eq:1.3}
\end{eqnarray}
Moreover, for $X\in \bar{L}^{0}$, we also have
\begin{equation}\label{eq:1.4}
E[-X|\cF_{t}]\leq -E[X|\cF_{t}].
\end{equation}
For the last inequality we used the convention $\infty-\infty=-\infty$.

Next, using \eqref{eq:1.1}-\eqref{eq:1.4}, we will prove the announced results. Assume that $X,Y\in\bar{L}^{0}$.

\noindent 1) If $\lambda\in L^{0}_{t}$, and $\lambda\geq 0$, then, by~(\ref{eq:1.1}) we get
\begin{align*}
E[\lambda  X|\cF_{t}] & = E[(\lambda X)^{+}|\cF_{t}]-E[(\lambda X)^{-}|\cF_{t}] =E[\lambda X^{+}|\cF_{t}]-E[\lambda X^{-}|\cF_{t}]=\\
& =\lambda  E[X^{+}|\cF_{t}]-\lambda E[X^{-}|\cF_{t}]= \lambda E[X|\cF_{t}].
\end{align*}
From here, and using~(\ref{eq:1.4}), for a general $\lambda\in L^{0}_{t}$,   we deduce
\begin{align*}
E[\lambda  X|\cF_{t}] & = E[1_{\{\lambda\geq 0\}}\lambda  X+1_{\{\lambda< 0\}}\lambda  X|\cF_{t}]= 1_{\{\lambda\geq 0\}}\lambda E[ X|\cF_{t}] +1_{\{\lambda< 0\}}(-\lambda) E[-X|\cF_{t}]\leq\\
& \leq 1_{\{\lambda\geq 0\}}\lambda E[ X|\cF_{t}] +1_{\{\lambda< 0\}}\lambda E[X|\cF_{t}]=\lambda E[X|\cF_{t}].
\end{align*}

\smallskip
\noindent 2) The proof of 2) follows from \eqref{eq:1.2} and \eqref{eq:1.4}; for $X\in L^{0}$ see also the proof in~\cite[Lemma 3.4]{Cherny2010}.

\smallskip
\noindent 3) On the set $\{E[X|\cF_{t}]=-\infty\} \cup \{E[Y|\cF_{t}]=-\infty\}$ the inequality is trivial due to the convention $\infty-\infty=-\infty$. On the other hand the set $\{E[X|\cF_{t}]>-\infty\} \cap \{E[Y|\cF_{t}]>-\infty\}$ can be represented as the union of the sets $\{E[X|\cF_{t}]>n\} \cap \{E[Y|\cF_{t}]>n\}$, for $n\in\bZ$, on which the inequality becomes the equality, due to \eqref{eq:1.3}.
\end{proof}

\end{appendix}

\section*{Acknowledgments}
Tomasz R. Bielecki and Igor Cialenco acknowledge support from the NSF grant  DMS-1211256.
Marcin Pitera acknowledges the support by Project operated within the Foundation for Polish Science IPP Programme ``Geometry and Topology in Physical Model'' co-financed by the EU European Regional Development Fund, Operational Program Innovative Economy 2007-2013. Part of the research was performed while Igor Cialenco was visiting the Institute for Pure and Applied Mathematics (IPAM), which is supported by the National Science Foundation.

The authors express their gratitude to the referees and the Associate Editor for valuable comments and suggestions that helped us to improve the paper.

{\small
\bibliographystyle{alpha}
\newcommand{\etalchar}[1]{$^{#1}$}

}

 \end{document}